\documentclass{amsart}

\usepackage{amssymb}
\usepackage{amsthm}
\usepackage{amsmath}
\usepackage{graphicx}
\usepackage{subfigure}
\usepackage{longtable}
\usepackage{psfrag}
\usepackage[noadjust]{cite}
\newtheorem{theorem}[subsection]{Theorem}
\theoremstyle{definition}
\newtheorem{definition}[subsection]{Definition}

\newcommand{\F}{\mathbb{F}}
\newcommand{\Z}{\mathbb{Z}}
\newcommand{\Q}{\mathbb{Q}}

\newcommand{\C}{\mathbb{C}}
\newcommand{\cH}{\mathcal{H}}

\DeclareMathOperator{\Frob}{Frob}
\DeclareMathOperator{\Gal}{Gal}

\newcommand{\GL}{\mathrm{GL}}
\newcommand{\SL}{\mathrm{SL}}

\newcommand{\Sym}{\mathrm{Sym}}
\newcommand{\aaa}{}

\begin{document}

\title[$\Sym^g$ cohomology for $\SL_{4}(\Z )$ and Galois representations]
{Cohomology with $\Sym^g$ coefficients for congruence subgroups of $\SL_4(\Z)$ and Galois representations}

\author{Avner Ash} \address{Boston College\\ Chestnut Hill, MA 02445}
\email{Avner.Ash@bc.edu} \author{Paul E. Gunnells}
\address{University of Massachusetts Amherst\\ Amherst, MA 01003}
\email{gunnells@umass.edu} \author{Mark McConnell}
\address{Princeton University\\ Princeton, New Jersey 08540}
\email{markwm@princeton.edu}

\thanks{
PG wishes to thank the National Science Foundation and the Simons Foundation}

\keywords{Cohomology of arithmetic groups, Galois representations, Voronoi complex, Steinberg
module, modular symbols}

\subjclass{Primary 11F75; Secondary 11F67, 20J06, 20E42}

\begin{abstract}
We extend the computations in \cite{AGM1, AGM2, AGM3,AGM4,AGM7} to find the
cohomology in degree five of a congruence subgroup~$\Gamma$ of
$\SL_4(\Z)$ with coefficients in~$\Sym^g(K^4)$, twisted by a nebentype
character~$\eta$, along with the action of the Hecke algebra.  This is
the top cuspidal degree. 
In this paper we take $K=\F$, a finite field of
large characteristic, as a proxy for~$\C$.  
For each Hecke eigenclass
found, we produce the unique Galois representation that appears to be
attached to it.  
%Our  computations show that in every case this Galois representation is the only one that could be attached to it.
%The existence of the attached Galois representations agrees with a theorem of Scholze \cite{scholze} (see also \cite{HLTT})
%and sheds light on the Borel-Serre boundary for~$\Gamma$.

The computations require modifications to our previous
algorithms to accommodate the fact that the coefficients are not one-dimensional.  %Nontrivialcoefficients add a layer of complication to our data structures, and
Types of attached Galois representations arise that were not found in our previous papers, and we must modify the Galois
Finder accordingly.
%We have also
%improved the Galois Finder so that it reports when the attached Galois
%representation is uniquely determined by our data.
\end{abstract}

\maketitle

\section{Introduction}\label{intro}

%\subsection{} 
The cohomology of arithmetic groups plays various roles in modern
number theory.  One of these concerns the connections between Hecke
eigenclasses in the cohomology and Galois representations.  This paper
continues our series of computations in this area for subgroups of
$\SL_4(\Z)$.

It is a highly nontrivial problem to compute the homology groups and
the action of the Hecke operators on them.  In this paper, we use the
sharbly complex for these computations, as we have done in the
previous papers of this series.

\begin{definition}
For any field $K$,
$\Sym^g(K^4)$ denotes the space of homogeneous polynomials of degree $g$ on $K^4$.  If $\eta$ is a nebentype character valued in $K^\times$, 
$\Sym^g(K^4)_\eta$ is defined to be $\Sym^g(K^4)\otimes K_\eta$.
\end{definition}

Choose a level $N$ and let $\Gamma=\Gamma_0(N)\subseteq \SL_4(\Z)$ and
let $K$ be a field.  The coefficient modules we study are
$\Sym^g(K^4)_\eta$ with a nebentype character~$\eta$.  Our earlier
papers only considered one-dimensional coefficient modules, i.e.,
$g=0$.

Since we are interested in automorphic representations, we would like
to set $K=\C$.  However, in order to avoid the inaccuracy of floating
point numbers in our huge linear algebra computations, we instead set
$K=\F$, with~$\F$ denoting a finite field~$\F_{p^r}$, where $p$ is
some prime with five decimal digits and $r\geqslant1$.  This large finite
field can be viewed as a proxy for $\C$.

We are thus interested in $H^*(\Gamma, \Sym^g(\F^4)_\eta)$ for various $N,g>0$,
and $\eta$.  For reasons explained in our earlier papers, we have only
investigated $*=5$, and that is also the only degree we study in this
paper.  We compute $H^5(\Gamma, \Sym^g(\F^4)_\eta)$ as a module for a
finite subset $\cH$ of the tame Hecke algebra (which is commutative),
and we diagonalize the $\cH$-action on this vector space.

%We assume that $\cH$ is generated by all the Hecke operators
%$T_{\ell,k}$ for all prime $\ell\leqslant\ell_0$.  Note: this is not true
%for N=18.

Let $z\in H^*(\Gamma, \Sym^g(\F^4)_\eta)$ be an $\cH$-eigenclass.

\begin{definition}
We say that a representation $\rho$ of the absolute Galois group of
$\Q$ is \emph{attached} to $z$ with respect to $\cH$ if for every Hecke
operator $T\in\cH$, if $T$ is supported at the prime $\ell$, then its
eigenvalue is that predicted by the equality of the characteristic
polynomial of $\rho(\Frob_\ell)$ and the Hecke polynomial of $z$ at
$\ell$.

We say that $\rho$ is \emph{entirely attached} to $z$ if for every Hecke
operator $T$ in the whole tame Hecke algebra, if $T$ is supported at
the prime $\ell$, then its eigenvalue is so predicted.
\end{definition}
  
In our previous papers we used the terminology that $\rho$ ``appears"
to be attached to $z$, because we can only compute a finite number of
Hecke operators $T$.  The new terminology contains more information by including the set of $T$ for which computations were made.  Also,
since Scholze has proved that there always exists a $\rho$ attached to
$z$, if there were a unique $\rho$ attached to $z$ with respect to $\cH$,
then $\rho$ must be entirely attached to $z$.  Of course, since we can find only a small number of Hecke eigenvalues (owing to the size of the matrices involved), even though our Galois finder returns only one $\rho$, it is still possible (though very unlikely) that this $\rho$ is an ``imposter'' and the entirely attached $\rho$ is some other Galois representation 
that agrees with $\rho$ on $\Frob_\ell$ for small $\ell$.

For each Hecke eigenclass $z$ computed in this paper, we find a
Galois representation attached to it, and this attached Galois
representation is uniquely determined by our data, in a sense to be
explained in Theorem~\ref{first}.  As we just explained, it is logically
possible (but very unlikely) that if we considered candidates for $\rho$ not in
the list of Section~\ref{basic} we might find other $\rho$'s
also attached to $z$ with respect to $\cH$.

If we could compute enough Hecke operators, then, using Scholze's
theorem and the method of Faltings--Serre, we could prove that a
given $\rho$ is entirely attached to a given $z$.  But it is not feasible to compute anywhere
near enough Hecke operators to do this for the homology classes found in this paper.

As in our earlier papers, these computations give new examples of
Scholze's theorem (recalled in~Section \ref{basic}) and new \textit{a posteriori}
tests of Serre-type conjectures for $\GL_4$.  
 
The Galois representations in this paper are all
reducible.  We do not know why certain combinations of characters and
cusp forms appear and others do not.  This ignorance stems from the
fact that the cohomology of $\overline X/\Gamma$ and its boundary, let
alone the restriction map from one to the other, is not known.  

Our computations are complete for the following values of
$(N,g,\eta)$.  For $g=1$ and~$2$, they are complete for all levels
$N\leqslant 18$, both prime and composite.  For $g=3,\dots,7$, as the
computations became slower, we computed only for certain prime levels
$N\leqslant 17$.  For $N=1$, we computed for $g\leqslant 10$.  When we
computed for a given pair $(N,g)$, we computed for all the~$\eta$
relevant to that pair.

The existence of attached Galois representations helps to corroborate
the correctness of our computations.  It is unimaginable that
attached Galois representations could be found if the computed Hecke
eigenvalues were random collections of numbers that had been calculated erroneously.

This paper covers the same ground as \cite{AGM7}, except for the difference
in the coefficient modules.  Therefore we refer the reader to \cite{AGM7} for
most of the background information and the description of how the
computations are performed.  We will explain below the changes needed in
order to deal with coefficient modules of dimension greater than 1,
and modifications required in the Galois Finder.  Then we will
summarize our findings and provide complete tables of the results of
our computations.

\section{Definitions, notations, basic constructions} \label{basic} 

\begin{definition} Fix $N\geqslant 1$.  

$\Gamma=\Gamma_0(N)$ will denote the
subgroup of matrices in $\SL_n(\Z)$ whose bottom row is congruent to
$(0,\dots,0,{*})$ modulo~$N$.

$\eta$ will denote a character of $(\Z/N)^\times$.  It can be viewed
as a character of $\Gamma$ by being applied to the $(n,n)$-entry of an
element in $\Gamma$.

$V$ will denote the standard representation of $\GL_4$.

$Sh_\bullet$ denotes the sharbly resolution of the Steinberg module for $\GL_4$.
\end{definition}

Recall that if $K$ is a field, and $\eta$ is a $K^\times$-valued character of
$(\Z/N)^\times$, we defined $K_{\eta}$ to be the one-di\-men\-sion\-al
vector space $K$ regarded as a $\Gamma$-module with action via the
nebentype character~$\eta$.  (We call $\eta$ the \emph{nebentype}
even if it is trivial.)

Section 2 of \cite{AGM7} gives the definitions of the Steinberg
module, the sharbly complex, and of the Hecke
polynomial at $\ell$.  It also explains why the sharbly homology is
isomorphic to $H^*(\Gamma, \Sym^g(\F^4)_\eta)$.  Section 3 of the same paper reviews
how the sharbly homology is calculated as a Hecke module, and Section
4 describes the Galois Finder. We will assume knowledge of these
matters in what follows.

There is an isomorphism of Hecke modules
\[
H^5(\Gamma,M)\approx H_1(\Gamma, Sh_{\bullet}\otimes_{\Z} M),
\]
where $M$ is any module on which the orders of the finite subgroups in
$\SL_{4} (\Z )$ are invertible; this condition is satisfied for us
since we will take $V=\F_{p^r}^4$ and $M =\Sym^g(V)\otimes \eta $,
where $p>5$.  

Indeed, in order to avoid the inaccuracy of floating point numbers in our huge
linear algebra computations, we use a finite field~$\F = \F_{p^r}$ as
a proxy for $\C$.  If $p>5$ and if there is no $p$-torsion in the
$\Z$-cohomology (which is very likely the case for large random $p$), then the $\C$- and mod~$p$-betti numbers coincide.
We use primes that have five decimal digits.
We choose $p$ and $r$ as follows.  

We
choose~$p$ so that the exponent of $(\Z/N)^\times$ divides $p-1$.
This makes the group of characters $(\Z/N)^\times\to\F_p^\times$
isomorphic to the group of characters $(\Z/N)^\times\to\C^\times$.  
(Note: this is not needed if $\eta=1$.  Some of our initial computations for 
$\eta=1$ were performed for a prime that differs from the $p$ we  used at the same level for nontrivial $\eta$.)
Later in the computation, we
choose~$r$ to ensure that the various Hecke eigenvalues that we compute lie in
$\F$. 

Define $S_{pN}$ to be the subsemigroup of integral matrices in
$\GL_n(\Q)$ satisfying the same congruence conditions mod~$N$ as $\Gamma$
and having positive determinant relatively prime to~$pN$.  Then
$\cH(pN)$, the \emph{tame Hecke algebra}, is the $\Z$-algebra of
double cosets $\Gamma S_{pN}\Gamma$.  It is a
commutative algebra that acts on the cohomology and homology of
$\Gamma$ with coefficients in any $S_{pN}$-module.  
$\cH(pN)$ is generated by all double cosets of the form
$\Gamma D(\ell,k)\Gamma$, where $\ell$ is a prime not
dividing $pN$, $0\leqslant k\leqslant n$, and $D(\ell,k)$
is the diagonal matrix with the first $n-k$ diagonal entries equal to
1 and the last $k$ diagonal entries equal to $\ell$.   When we
consider the double coset generated by $D(\ell,k)$ as a Hecke
operator, we call it $T(\ell,k)$.

${\F_\eta}$ is an $S_{pN}$-module, where a matrix $s\in S_{pN}$
acts on $\F$ via $\eta(s_{nn})$, where~$s_{nn}$ is the lower right entry of the $n\times n$ matrix $s$.

\begin{definition}\label{def:hp}
Let $V$ be an $\F[\cH(pN)]$-module.  Suppose that $v\in V$ is a
simultaneous eigenvector for all $T(\ell,k)$ and that
$T(\ell,k)v=a(\ell,k)v$ with $a(\ell,k)\in\F$ for all prime $\ell
\nmid pN$ and $0\leqslant k\leqslant n$.  If
\[
\rho\colon G_\Q\to \GL_n(\F)
\]
is a continuous representation of
$G_{\Q} = \Gal (\overline\Q/\Q)$ unramified outside $pN$, and if
\begin{equation}\label{eqn:hp}
\sum_{k=0}^{n}(-1)^k\ell^{k(k-1)/2}a(\ell,k)X^k=\det(I-\rho(\Frob_\ell)X)
\end{equation}
for all $\ell\nmid \ pN$, then we say that~$\rho$ is \emph{entirely attached} to~$v$. 
\end{definition}  
Here, $\Frob_\ell$ refers to an arithmetic Frobenius element, so that
if $\varepsilon$ is the cyclotomic character, we have $\varepsilon
(\Frob_\ell)=\ell$.  

The polynomial on the left-hand side of~(\ref{eqn:hp}) is called the \emph{Hecke polynomial}
for~$v$ at~$\ell$.

\begin{definition} 
Suppose we have $v$ and $\rho$ as above.  If $\cH$ is a subset of 
$\cH(pN)$ and the eigenvalues of all $T\in\cH$ are those that would be predicted by $\rho$ if $\rho$ were entirely attached to $v$, then we say 
$\rho$ is \emph{attached} to~$v$ with respect to $\cH$. 
\end{definition} 

The following  is a special case
of a theorem of Scholze:
\begin{theorem}\label{scholze} 
Let $N\geqslant1$.  
Let $v$ be a Hecke eigenclass in
$H^{5}(\Gamma_0(N), \Sym^g(\F^4)_\eta)$.  Then
there is entirely attached to~$v$ a continuous representation $\rho$, unramified
outside $pN$:
\[
\rho\colon G_\Q\to \GL_n( \F).
\]
\end{theorem}
\noindent Since $\rho$ is entirely attached to $v$, it is unique up to isomorphism.

The coefficient modules $M$ studied in this paper are 
$\Sym^g(\F^4)_\eta$ for various levels $N$, nebentypes $\eta$ and degrees $g$.
We compute homology and the Hecke action exactly as in \cite{AGM7}.  Of
course we have to modify the programs to use $\Sym^g(\F^4)_\eta$
coefficients.

When we wrote our code for \cite{AGM7}, we had made sure to support
arbitrary coefficient modules $M$ for the cohomology.  During
\cite{AGM7}, this code was tested for $M$'s that were one-dimensional over
$\mathbb{F}$.  As it turned out, the same code worked out of the box
for the high-dimensional $M$ used in the present paper, after 
a few small incompatibility bugs were fixed.

As explained in \cite{AGM3,AGM7}, computing the cohomology comes down
to finding the kernels and images of certain large matrices coming
from the cells of the well-rounded retract.  The dimension of
$\Sym^g(\mathbb{F}^4)$ is $\binom{g+3}{3} \sim g^3/6$.  Thus, when $M
= \Sym^g(\mathbb{F}^4) \otimes \eta$, the numbers of rows and columns
in the matrices are larger by approximately this factor of $g^3/6$,
compared to the size for $M=\F_\eta$ alone.  In turn, the size for
one-dimensional~$M$ grows roughly like $O(N^3)$ for both rows and
columns.  This explains why we stopped our computations at $g=7$, and
why we restricted ourselves to smaller ranges of $g$ as $N$ became
large (or a slightly larger range for $N=1$).  The largest matrix we
encountered was for $N=18$ and $\Sym^2(V)$, where the matrix was
$16204\times 56420$.  This is far smaller than the largest matrix in
\cite{AGM3}, which was about 1~million by 4~million (for the case
$N=211$, $M = 1$).  However, in \cite{AGM3}, unlike in \cite{AGM7} and
the present paper, we were not computing bases for the kernels and
images of the matrices, which are needed in order to compute the Hecke
operators; we were only computing ranks of matrices mod~$p$.  That is
why we could go to much larger matrices in \cite{AGM3}.

To find attached Galois representations, we use the Galois Finder
program, which is part of our Sage code.  We had to modify it for the current project, making  two changes.  First, it now
considers cusp forms of all weights $2, \dots, g+4$.  In \cite{AGM7}, where
$g=0$, we only needed to consider weights $2,3,4$.  Secondly, it
now considers powers $\varepsilon^i$ of~$\varepsilon$ for all~$i = 0,
\dots, g+3$, as opposed to $0,\dots,3$ for \cite{AGM7}.

We compute 
the action 
on $V=H_1(\Gamma_0(N), Sh_{\bullet} \otimes_\Z \Sym^g(\F^4)_\eta)$
of the Hecke operators $T(\ell,k)$ for $k=1, 2, 3$ and for~$\ell$
ranging through a set
\[
L = \bigl\{\ell \bigm | \ell \mathrm{\ prime,} \, \ell \leqslant \ell_0, \,
\ell\nmid pN\bigr\}.
\]
The upper bound~$\ell_0$ depends on the level~$N$ and the
nebentype~$\eta$.  What limits the choice of $\ell_0$ is the size of
the matrices involved in the computation and the time it takes.

In this paper,  $5\leqslant\ell_0\leqslant11$.  For $\ell_0$ itself we sometimes find only $T(\ell_0,1)$ and not $T(\ell_0,k)$
for $k=2,3$ because of the size of the computations.   For $k=0,4$, we do not have to do any computation: $T(\ell,0)$ is the identity and $T(\ell,4)$ is
$\eta(\ell)\ell^g$ times the identity.  To check our work, we always verify
that the Hecke operators commute pairwise.

\section{Observations from the data}

In the range of our computations, all the Galois representations that occur are reducible with constituents of dimension 1 and 2.
One-\-di\-men\-sion\-al constituents come from Dirichlet characters mod~$N$ taking
values in the cyclotomic field $K_0$ generated by a primitive $N$-th root of unity.  
Two-\-di\-men\-sion\-al constituents come  from newforms of level dividing~$N$ and
weights $2,\dotsc ,g+4$.  
Any of these constituents may be multiplied by a power of the cyclotomic character.

Let $K_1, K_2, \dots$  be the fields of coefficients of the $q$-expansions of the newforms we
have listed, together with~$K_0$.  The Galois Finder works in the residue class fields for the various
primes $\mathfrak{P}$ over~$p$ in the various~$K_i$'s.  We define~$r$
to be the smallest integer so that all these residue class fields
embed in $\F = \F_{p^r}$.  We choose $p$ to make $r$ as small as possible, given the constraint that $p$ should be no more than five digits (which is needed for speed).  The field~$\F$ is recorded at the top of
each table in Section~\ref{sec:results}.
The table also specifies, for each $N,\eta,g$, the set of Hecke operators making up our choice of $\cH$ for those parameters.

We summarize our first observation as follows:

\begin{theorem}\label{first}
For $N$, $p^r$, $\eta$, $\ell$, and $g$ as covered in the tables in Section~\ref{sec:results}, 
the Hecke operators $T\in\cH$ on $H^5(\Gamma_0(N), {\Sym^g(\F_{p^r})}_\eta)$ are all semisimple.  For every Hecke eigenvector $z$, there exists a unique reducible Galois representation $\rho:G_\Q\to \GL_4(\F_{p^r})$ (within the scope of the Galois finder) that is attached to $z$ with respect to $\cH$. Each such $\rho$ is either the sum of four characters or the sum of two characters plus the Galois representation of a newform tensored with a character.
\end{theorem}

Let $E$ denote a simultaneous eigenspace of $\cH$ on $V
=H^5(\Gamma_0(N), \Sym^g(\F^4)_\eta)$, where $\F= \F_{p^r}$.
We define two kinds of multiplicity for~$E$.  
\begin{definition}
The
\emph{Hecke multiplicity} of~$E$ equals $\dim_{\F} E$.
\end{definition}
Let $G_\eta$ be the stabilizer of $\eta$ in the Galois group of  
$\F/\F_p$.  Then  $G_\eta$ acts on $V$ and permutes the Hecke eigenspaces.
\begin{definition}
The
\emph{Galois multiplicity} of~$E$ equals the cardinality of the orbit of $E$ under $G_\eta$. 
\end{definition}

The Galois finder works exactly as it did in \cite{AGM7}.  If the
Galois finder returns the same $\rho$ exactly $d$ times, for Hecke
eigenspaces $E_1,\dots,E_d$, then the Galois multiplicity of each
$E_i$ equals $d$ and we list only one of them in the tables.  Although
it seems like it is returning the same $\rho$, this is not true: it is
using a different prime $\frak P$ for each one.

We now describe in detail the list of Galois
representations~$\rho$ which our Galois Finder used for this paper.

First are the Dirichlet characters~$\chi$ with values in~$\F$, which
we identify with one-di\-men\-sion\-al Galois representations as
usual.  We take all the characters of conductor $N_1$ for all $N_1
\mid N$.  Sage's class \texttt{DirichletGroup} enumerates the~$\chi$
automatically.  The characteristic polynomial of Frobenius at~$\ell$
for~$\chi$ is $1 + \chi(\ell) X$, for all $\ell \nmid pN$.
%Each~$\chi$ can be lifted to characteristic zero, since $p\equiv 1
%\pmod{N}$.

Another one-di\-men\-sion\-al character is the cyclotomic
character~$\varepsilon$.  We look at $\varepsilon^w$ for $w=0,1,2\dots,g+3$,
because these are the powers predicted by the generalizations of Serre's conjecture
for mod~$p$ Galois representations \cite{AS,ADP}.  

We define the \emph{Hodge-Tate (HT) numbers} for Galois representations as follows.
For a character $\chi\otimes \varepsilon^w$, there is
a list of one integer $[w]$.  To a representation coming
from a newform~$\rho$ of weight~$k$, there is a list of two
integers, $[0, k-1]$.  More generally, for $\chi\otimes
\varepsilon^w\otimes\rho$, the list is $[w, w+k-1]$.  For direct sums
of representations, the lists are concatenated and then ordered by increasing values of the entries.  For the
four-di\-men\-sion\-al Galois representations we find that fit our data,
we always observe that the list is $[0,1,2,g+3]$ after sorting. 
This is
predicted by the Serre-type conjectures
% and the
%conjectural HT numbers.  
and gives us a check on our computations.
%See also Section \ref{HT}.

Another check on our computations comes from considering the
relationship between the nebentype character and the determinant of
the attached representation.  Suppose a Galois representation $\rho$
is attached to a Hecke eigenclass in $H^5(\Gamma_0(N),
\Sym^g(\F^4)_\eta)$.  Then the determinant of $\rho(\Frob_\ell)$ must
equal the coefficient of $X^4$ in the Hecke polynomial, namely
$\eta(\ell)\ell^{g+6}$.  We observe that this is always the case in
our data.  

\section{Other observed regularities in the data}\label{ChiReasons}

In this section, we set $V=\F^4$.
%we let $\Gamma_0(a,b)$ denote the subgroup of
%$\GL_a(\Z)$ where the bottom row is congruent to $(0,\dots,0,*)$
%modulo $b$.  Thus $\Gamma_0(N)=\Gamma_0(4,N)\cap \SL_4(\Z)$ in our
%notation.  We 
A Hecke eigenclass in
$H^5(\Gamma_0(N),\Sym^g(V)_\eta)$ will be denoted by the letter $z$, and its attached
Galois representation by~$\rho$.

\subsection{Oddness}\label{odd}
We observe that~$\rho$ is always odd. In other words,
the eigenvalues of $\rho(c)$ are $+1,-1,+1,-1$, where~$c$ denotes
complex conjugation.  This must be the case, as follows from a theorem
of Caraiani and LeHung \cite{CLH}.

\subsection{Multiplicities}\label{mult}

We observe that the Galois multiplicity of eigenspaces in our data can
be any integer from 1 to 6, while the Hecke multiplicity of
eigenspaces in our data can be 1, 3, 4, 6, or~9.  We do not have an explanation for why
other multiplicities do not occur.  It is possible that more
computations would reveal other multiplicities.

\subsection{Patterns}\label{patterns} Recall that $N$
denotes the level of the Hecke eigen\-space and $\eta$ denotes the
nebentype of the coefficients.

Each Galois representation in the tables is one of the following
types.  We let $\chi$ and $\psi$ denote 1-dimensional Galois
representations with conductor dividing $N$ and $\sigma_k$ an
irreducible 2-dimensional Galois representation corresponding to a
newform of weight $k$ and level dividing $N$.

1) $\chi\varepsilon^0\oplus \varepsilon^1 \oplus \varepsilon^2 \oplus \varepsilon^{g+3}$ and
$\varepsilon^0\oplus \varepsilon^1 \oplus \chi\varepsilon^2 \oplus \varepsilon^{g+3}$, $\chi\ne1$.   These always occur in such pairs.
  
  2)  $\chi\varepsilon^0\oplus \varepsilon^1 \oplus \varepsilon^2 \oplus \psi\varepsilon^{g+3}$ and
  $\varepsilon^0\oplus \varepsilon^1 \oplus \chi\varepsilon^2 \oplus \psi\varepsilon^{g+3}$, $\chi\ne1$, $\psi\ne1$.   These always occur in such pairs.  They are much rarer than type 1).

3) $\varepsilon^0\oplus \varepsilon^1 \oplus \varepsilon^2\sigma_{g+2}$.  
 (Notice that no summand here gets multiplied by a nontrivial character.)
 This always occurs unless there is no 
$\sigma_{g+2}$ with nebentype equal to $\eta$.  

4) $\varepsilon^1\oplus \varepsilon^{g+3} \oplus \varepsilon^0\sigma_{3}$.
(Notice that no summand here gets multiplied by a nontrivial character.)
Whenever type 4) occurs for given $N,g,\eta$, there also occurs type 3) and  type 1).

5) $\chi\varepsilon^0\oplus \varepsilon^{2} \oplus \varepsilon^1\sigma_{g+3}$ and
$\varepsilon^0\oplus \chi\varepsilon^{2} \oplus \varepsilon^1\sigma_{g+3}$,
$\chi\ne1$.   These always occur in such pairs.  This type occurs in our data for $N=12,15,16,18$.

\subsection{Differences from our previous findings for $g=0$}\label{diff}

Unlike in \cite{AGM7}, neither $\varepsilon^1$ nor $\sigma_k$ is ever
multiplied by a nontrivial character.  Of course, more data might
disturb this observation.

In \cite{AGM7}, if $\eta$ factors nontrivially as $\eta=\psi\chi$ then
either all three of the following or none of the following occur:

$\rho=  \psi\varepsilon^0\oplus \chi \varepsilon^1\oplus\varepsilon^2\oplus \varepsilon^3$

$\rho= \psi \varepsilon^0\oplus \varepsilon^1\oplus\varepsilon^2\oplus \chi\varepsilon^3$

$\rho=  \varepsilon^0\oplus \varepsilon^1\oplus\psi\varepsilon^2\oplus \chi\varepsilon^3$.

\noindent The natural analogue of this assertion for $g>0$ is  not true. See for example $N=16$.  

In \cite{AGM7}, 
 $\varepsilon^2 \chi\sigma$ ($\chi$ possibly trivial) occurred as a summand for a given coefficient module if and only if
 $\varepsilon^0
\chi\sigma$ occurred.  This is no longer true when $g>0$.  For example, see $N=3$, Coeffs $=\Sym^6(V)$.

\subsection{Heuristics}\label{h}

We do not have explanations for most of the regularities observed
above.  In \cite{AGM7} we gave a heuristic for the conductors of the
characters and the levels and weights of the cuspforms that appear in
the tables by referring to the homology of various parabolic subgroups
of $\GL_4$ intersected with $\Gamma$.  We refer to \cite{AGM7} for the
details of this analysis, and very briefly discuss them as they apply
to the tables below.  In \cite[Section 5.6]{AGM7} the analysis was
accompanied by five diagrams lettered (a) through (e), which we have
reproduced in Figure \ref{fig:specseq} below.

The heuristic concerning the conductors of the characters and the levels  of the cuspforms is the same as before.  As for the weights of the cuspforms:

The Borel--Serre boundary $B_\Gamma$ is the union of faces $F(P)$,
where $P$ runs over a set of representatives of $\Gamma$-orbits of
parabolic subgroups $P$ of $\GL_4(\Q)$.  Each parabolic subgroup~$P$
is conjugate to a standard one with block sizes $(n_1,\dots,n_{k+2})$
down the diagonal.  We call this tuple the ``type" of $P$ or of
$F(P)$.  The nonnegative integer $k$ equals the codimension of $F(P)$
in $B_\Gamma$.

Our heuristic explanation assumes that each Hecke eigenspace restricts
nontrivially to at least one of the faces.  Our data all conforms to assuming
this face is type $(2,2)$.  By the Eichler--Shimura theorem, a block of
size 2 will give rise to the Galois representation $\sigma$ attached
to a holomorphic cuspform with level dividing $N$, or to a sum of two
characters (in the case of an Eisenstein series), with conductors
dividing $N$.  In general, $\sigma$ and these characters may be
multiplied by a power of the cyclotomic character.

We now use this heuristic to describe the various kinds of Galois
representations that occur in our data in the tables.  Write the
parabolic subgroup of type $(2,2)$ as $P=L_1L_2U$ where $L_i\approx
\GL_2$ for $i=1,2$ and $U$ is the unipotent radical of $P$.  Note that
$\Sym^g$ restricted to $L_2$ has a submodule isomorphic to $\Sym^g$
for $\GL_2$, and another submodule isomorphic to $\Sym^0$ for $\GL_2$.

\begin{figure}[htb]
\psfrag{1}{$ 1$}
\psfrag{2}{$ 2$}
\psfrag{3}{$ 3$}
\psfrag{4}{$ 4$}
\psfrag{cusp}{$ \text{cusp}$}
\psfrag{eis}{$ \text{eis}$}
\centering
\subfigure[\label{fig:1a}]{\includegraphics[scale=0.22]{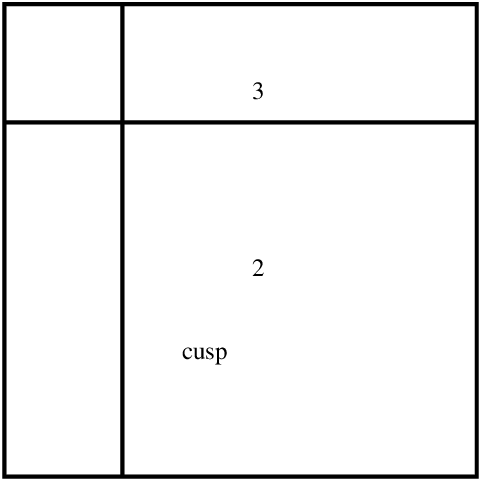}}
\quad\quad
\subfigure[\label{fig:1b}]{\includegraphics[scale=0.22]{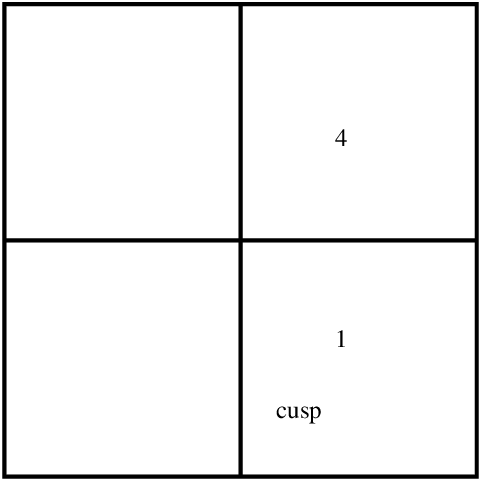}}
\quad\quad 
\subfigure[\label{fig:1c}]{\includegraphics[scale=0.22]{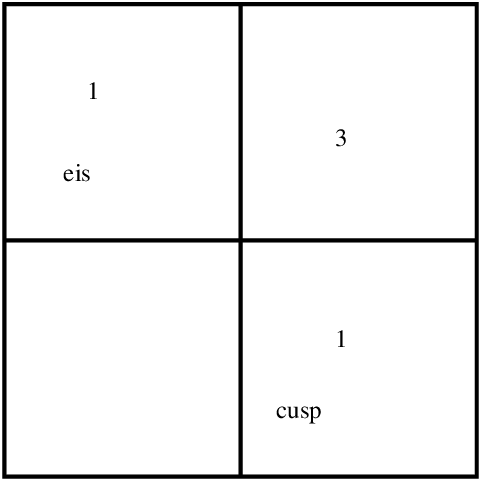}}
\quad\quad
\subfigure[\label{fig:1d}]{\includegraphics[scale=0.22]{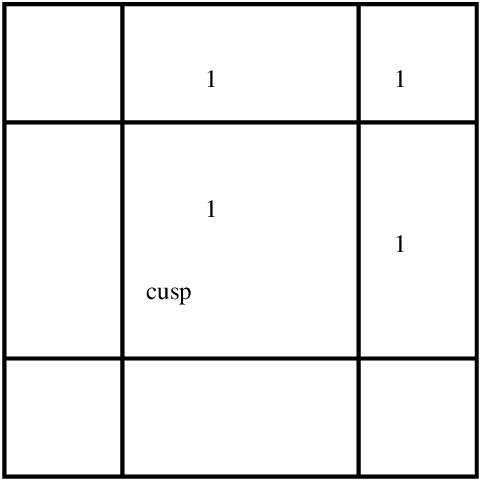}}
\quad\quad 
\subfigure[\label{fig:1e}]{\includegraphics[scale=0.22]{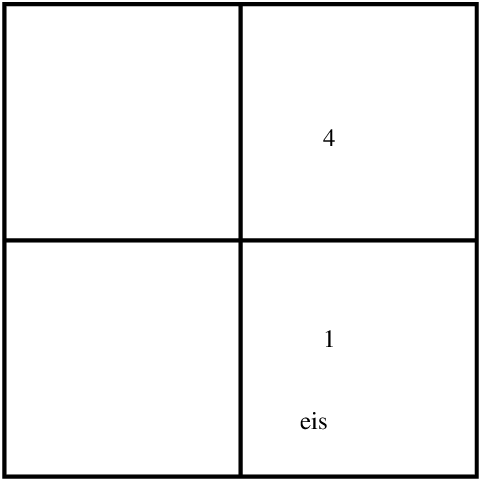}}
\caption{\label{fig:specseq} Schematics of homology classes on faces
 of the Borel--Serre boundary}
\end{figure}

\subsection{Holomorphic cusp forms of weight $g+2$}\label{wg+2}

In this case (Figure~1(b)), when we restrict the coefficients
$\Sym^g(V)\otimes\eta$ to $L_2$ we use the submodule $W_g$ isomorphic
to $\Sym^g$ for $\GL_2$.  We place a cuspform on the $L_2$ block of
weight $g+2$ (corresponding to the homology of the arithmetic group in
the $\GL_2$-block with coefficients in $H_4(U)\otimes W_g$.)  This
gives classes of type 3).

%We could place an Eisenstein series on the $L_1$ block which would
%lead to nontrivial characters appearing in the Galois representation.
%We do not know why this doesn't happen in our data.  Similar remarks
%apply below to other seemingly feasible patterns that do not occur.

\subsection{Holomorphic cusp forms of weight $g+3$}\label{wg+3}

In this case (Figure~1(c)), when we restrict the coefficients
$\Sym^g(V)\otimes\eta$ to $L_2$ we use the submodule $W_g$ isomorphic
to $\Sym^g$ for $GL_2$.  We place a cuspform on the $L_2$ block of
weight $g+3$ (corresponding to the homology of the arithmetic group in
the $GL_2$-block with coefficients in $H_3(U)\otimes W_g$.)  We place
an Eisenstein series on the $L_1$ block. This gives classes of type
5).

\subsection{Holomorphic cusp forms of weight 3}\label{w3} In this case
(also Figure~1(c)), when we restrict the coefficients
$\Sym^g(V)\otimes\eta$ to $L_2$ we use the submodule $W_0$ isomorphic
to $\Sym^0$ for $GL_2$.  We place a cuspform on the $L_2$ block of
weight $3$ (corresponding to the homology of the arithmetic group in
the $GL_2$-block with coefficients in $H_3(U)\otimes W_0$.)  
We place
an Eisenstein series on the $L_1$ block.This
gives classes of type 4).

\subsection{Sums of 4 characters}\label{allchars}

In this case (Figure~1(e)), when we restrict the coefficients
$\Sym^g(V)\otimes\eta$ to $L_2$ we use the submodule $W_g$ isomorphic
to $\Sym^g$ for $GL_2$.  We place an Eisenstein series on the $L_2$
block of weight $g+3$ (corresponding to the homology of the arithmetic
group in the $GL_2$-block with coefficients in $H_3(U)\otimes W_g$.)
%We also place an Eisenstein series on the $L_1$ block.
This gives classes of types 1) and 2).  

We do not understand the finer
details of the sums of characters, nor of the other types of attached Galois representations.

\subsection{Missing patterns from \cite{AGM7}}\label{miss} Nothing in
our data corresponds to Figure~1(d) where $P$ is a
$(1,2,1)$-parabolic subgroup nor to Figure~1(a) where $P$ is a
$(1,3)$-parabolic subgroup.  We would expect Figure~1(a) to occur if $\GL_3$ has a cuspidal cohomology class of level $N$, but such $N$ are beyond the range of our computations.
We do not have a guess as to whether Figure~1(d) would occur for larger
levels $N$.

\section{Tables of results}\label{sec:results}

\subsection{} \label{resMain}
The tables in this section present the main results of the paper.

Let $V$ denote the {standard representation} of $\GL_4$ acting on a
vector space of dimension 4.  A given coefficient module will be
denoted $\Sym^g(V)\otimes\eta$ for a nebentype~$\eta$.  (We used a
subscript $\eta$ earlier, but putting $\eta$ on the line makes it
easier to read).  Dirichlet characters will be denoted by a
subscripted $\chi$, and 2-dimensional irreducible Galois
representations will be denoted by a subscripted $\sigma$.  We replace
$\Sym^1(V)$ with~$V$ in the tables.

The topmost box in each table gives the level~$N$, the coefficient
module $\Sym^k(V)\otimes\eta$ with nebentype~$\eta$, and the
field~$\F_{p^r} = GF(p^r)$ that was our proxy for~$\C$.  We include
only one representative for each Galois orbit of nebentype characters.
Next, we list the Hecke operators we computed.  $T_\ell$ means we
computed $T_{\ell,1}$, $T_{\ell,2}$, and $T_{\ell,3}$.  Listing
$T_{\ell,1}$ means we computed only that part of $T_\ell$.

The succeeding rows in each table give the Galois multiplicity
(Def.~2.8), the Hecke multiplicity (Def.~2.7), and the Galois
representation itself.
%% that is found by the Galois Finder to be apparently attached to the
%% given Hecke eigenclass in the cohomology
%% $H^5(\Gamma_0(N),\F_\eta)$.
The cyclotomic character is denoted~$\varepsilon$.
  
$\chi_{N}$ or~$\chi_{N,i}$ are a basis for the mod~$p$ Dirichlet
characters $(\Z/N\Z)^\times \to \F_p$.  They are listed in a separate
table at the end.  As explained above, we
usually\footnote{In some computations with trivial nebentype $\eta=1$,
we were not concerned with the Dirichlet characters as a group.  In these cases, we
arbitrarily chose $p=12379$, the fourth prime after $12345$.}
choose~$p$ depending on~$N$ so that the exponent of $(\Z/N\Z)^\times$
divides the order $p-1$ of $\F_p^\times$.  It follows that the group
of complex-valued Dirichlet characters is isomorphic to the group of
mod~$p$ characters.

The symbol $\sigma_{N.k.\mathrm{a}.\mathrm{x}}$ denotes a classical
cuspidal holomorphic newform.  We label these following the
conventions of the LMFDB~\cite{lmfdb}.
%%% bibtex reference to lmfdb is in comments below
Thus $N$ is the level of the newform, $k$ is its weight,
``$\mathrm{a}$'' is the LMFDB name for the nebentype character of the
newform, and ``$\mathrm{x}$'' denotes a specified Galois orbit of
newforms.  We use the same symbol $\sigma_{N.k.\mathrm{a}.\mathrm{x}}$
to stand for the two-di\-men\-sional Galois representation attached to
the cusp form of that name.

%%%%%%%%%%%%%%%%%%%%%%%%%%%%%%%%%%%%%%%%%%%%%%%%%%%%%%%%%%%

%% @misc{lmfdb,
%%   shorthand    = {LMFDB},
%%   author       = {The {LMFDB Collaboration}},
%%   title        = {The {L}-functions and modular forms database},
%%   howpublished = {\url{https://www.lmfdb.org}},
%%   year         = {2023},
%%   note         = {[Online; accessed 30 July 2023]},
%% }$N$ is the level, $k$ is the weight, and $N.a$ is the

For $g=1$ and~$2$, we computed the cohomology for all levels
$N\leqslant 18$, both prime and composite.  For
$g=3,\dots,7$, as the computations became slower, we computed only for
certain prime levels $N\leqslant 17$.  For $N=1$, we computed for
$g\leqslant 10$.  When we computed for a given pair $(N,g)$, we
computed for all the~$\eta$ relevant to that pair.  In general, the
range of~$N$ for which we computed became smaller as~$g$ grew larger.

\subsubsection*{Level 1}

\begin{center}
\begin{tabular}{|l|}
  \hline
  \textbf{Level $N = 1$.  Coeffs $\Sym^g(V)$ for $g=2,4,6,8$.  Field $\mathbb{F} = GF(12379)$.} \\
  \hline
  Dim~$0$. \\
  \hline
  \end{tabular}
\end{center}

\begin{center}
\begin{tabular}{|r|r|c|}
\hline
\multicolumn{3}{|l|}{\textbf{Level $N = 1$.  Coeffs $\Sym^{10}(V)$.  Field $\mathbb{F} = GF(12379)$.}} \\
\multicolumn{3}{|l|}{Computed $T_{2}$, $T_{3}$.  Dim~$1$.} \\ % del-longform
% \hline %% del-shortform
%  & $T_{2}$ & $(1-x)(1-2x)(1-3313x)(1+3409x)$ \\ % del-shortform
%  & $T_{3}$ & $(1-x)(1-3x)(1-2268x+1646x^{2})$ \\ % del-shortform
\hline
1 & $1$ & $\aaa \varepsilon^{0} \oplus \aaa \varepsilon^{1} \oplus \aaa \varepsilon^{2} \sigma_{1.12.\mathrm{a}.\mathrm{a}}$ \\
%% \hline
%% \multicolumn{3}{|l|}{$\sigma_{1} = \mathrm{Cusp\ form\ }q - 24q^{2} + 252q^{3} - 1472q^{4} + 4830q^{5} - 6048q^{6} - 16744q^{7} + O(q^{8})$} \\
%% \multicolumn{3}{|l|}{$\mathrm{for\ }\mbox{\rm SL}_2(\Bold{Z})\mathrm{\ with\ character\ }\aaa\mathrm{\ in\ }\mathbb{F}$} \\
%% \multicolumn{3}{|l|}{$\mathrm{of\ weight\ }12\mathrm{\ with\ coeffs\ in\ }\mathbb{Q}$} \\
\hline
\end{tabular}
\end{center}

\subsubsection*{Level 2}

\begin{center}
  \begin{tabular}{|l|}
  \hline
  \textbf{Level $N = 2$.  Coeffs $\Sym^2(V)$.  Field $\mathbb{F} = GF(12379)$.} \\
  \hline
  Dim~$0$. \\
  \hline
  \end{tabular}
\end{center}

\begin{center}
  \begin{tabular}{|l|}
  \hline
  \textbf{Level $N = 2$.  Coeffs $\Sym^4(V)$.  Field $\mathbb{F} = GF(12379)$.} \\
  \hline
  Dim~$0$. \\
  \hline
  \end{tabular}
\end{center}

\begin{center}
\begin{tabular}{|r|r|c|}
\hline
\multicolumn{3}{|l|}{\textbf{Level $N = 2$.  Coeffs $\Sym^6(V)$.  Field $\mathbb{F} = GF(12379)$.}} \\
\multicolumn{3}{|l|}{Computed $T_{3}$, $T_{5}$, $T_{7,1}$.  Dim~$1$.} \\ % del-longform
% \hline %% del-shortform
%  & $T_{3}$ & $(1-x)(1-3x)(1-108x+3841x^{2})$ \\ % del-shortform
%  & $T_{5}$ & $(1-x)(1-5x)(1-2011x)(1-5118x)$ \\ % del-shortform
%  & $T_{7}$ & $1-276x + O(x^2)$ \\ % del-shortform
\hline
1 & $1$ & $\aaa \varepsilon^{0} \oplus \aaa \varepsilon^{1} \oplus \aaa \varepsilon^{2} \sigma_{2.8.\mathrm{a}.\mathrm{a}}$ \\
%% \hline
%% \multicolumn{3}{|l|}{$\sigma_{2} = \mathrm{Cusp\ form\ }q - 8q^{2} + 12q^{3} + 64q^{4} - 210q^{5} - 96q^{6} + 1016q^{7} + O(q^{8})$} \\
%% \multicolumn{3}{|l|}{$\mathrm{for\ }\Gamma_0(2)\mathrm{\ with\ character\ }\aaa\mathrm{\ in\ }\mathbb{F}$} \\
%% \multicolumn{3}{|l|}{$\mathrm{of\ weight\ }8\mathrm{\ with\ coeffs\ in\ }\mathbb{Q}$} \\
\hline
\end{tabular}
\end{center}

\subsubsection*{Level 3}

\begin{center}
\begin{tabular}{|r|r|c|}
\hline
\multicolumn{3}{|l|}{\textbf{Level $N = 3$.  Coeffs $V\otimes\chi_3$.  Field $\mathbb{F} = GF(12379)$.}} \\
\multicolumn{3}{|l|}{Computed $T_{2}$, $T_{5}$, $T_{7}$, $T_{11}$.  Dim~$2$.} \\ % del-longform
% \hline %% del-shortform
%  & $T_{2}$ & $(1-x)(1-2x)(1+4x)(1-16x)$ \\ % del-shortform
%  & $T_{5}$ & $(1-x)(1-5x)(1+25x)(1-625x)$ \\ % del-shortform
%  & $T_{7}$ & $(1-x)(1-7x)(1-49x)(1-2401x)$ \\ % del-shortform
%  & $T_{11}$ & $(1-x)(1-11x)(1+121x)(1-2262x)$ \\ % del-shortform
\hline
1 & $1$ & $\aaa \varepsilon^{0} \oplus \aaa \varepsilon^{1} \oplus \chi_{3} \varepsilon^{2} \oplus \aaa \varepsilon^{4}$ \\
% \hline %% del-shortform
% \hline %% del-shortform
%  & $T_{2}$ & $(1+x)(1-2x)(1-4x)(1-16x)$ \\ % del-shortform
%  & $T_{5}$ & $(1+x)(1-5x)(1-25x)(1-625x)$ \\ % del-shortform
%  & $T_{7}$ & $(1-x)(1-7x)(1-49x)(1-2401x)$ \\ % del-shortform
%  & $T_{11}$ & $(1+x)(1-11x)(1-121x)(1-2262x)$ \\ % del-shortform
\hline
1 & $1$ & $\chi_{3} \varepsilon^{0} \oplus \aaa \varepsilon^{1} \oplus \aaa \varepsilon^{2} \oplus \aaa \varepsilon^{4}$ \\
\hline
\end{tabular}
\end{center}

\begin{center}
  \begin{tabular}{|l|}
  \hline
  \textbf{Level $N = 3$.  Coeffs $\Sym^2(V)$.  Field $\mathbb{F} = GF(12379)$.} \\
  \hline
  Dim~$0$. \\
  \hline
  \end{tabular}
\end{center}

\begin{center}
\begin{tabular}{|r|r|c|}
\hline
\multicolumn{3}{|l|}{\textbf{Level $N = 3$.  Coeffs $\Sym^3(V)\otimes\chi_3$.  Field $\mathbb{F} = GF(12379)$.}} \\
\multicolumn{3}{|l|}{Computed $T_{2}$, $T_{5}$, $T_{7}$.  Dim~$2$.} \\ % del-longform
% \hline %% del-shortform
%  & $T_{2}$ & $(1-x)(1-2x)(1+4x)(1-64x)$ \\ % del-shortform
%  & $T_{5}$ & $(1-x)(1-5x)(1+25x)(1-3246x)$ \\ % del-shortform
%  & $T_{7}$ & $(1-x)(1-7x)(1-49x)(1+6141x)$ \\ % del-shortform
\hline
1 & $1$ & $\aaa \varepsilon^{0} \oplus \aaa \varepsilon^{1} \oplus \chi_{3} \varepsilon^{2} \oplus \aaa \varepsilon^{6}$ \\
% \hline %% del-shortform
% \hline %% del-shortform
%  & $T_{2}$ & $(1+x)(1-2x)(1-4x)(1-64x)$ \\ % del-shortform
%  & $T_{5}$ & $(1+x)(1-5x)(1-25x)(1-3246x)$ \\ % del-shortform
%  & $T_{7}$ & $(1-x)(1-7x)(1-49x)(1+6141x)$ \\ % del-shortform
\hline
1 & $1$ & $\chi_{3} \varepsilon^{0} \oplus \aaa \varepsilon^{1} \oplus \aaa \varepsilon^{2} \oplus \aaa \varepsilon^{6}$ \\
\hline
\end{tabular}
\end{center}

\begin{center}
\begin{tabular}{|r|r|c|}
\hline
\multicolumn{3}{|l|}{\textbf{Level $N = 3$.  Coeffs $\Sym^4(V)$.  Field $\mathbb{F} = GF(12379)$.}} \\
\multicolumn{3}{|l|}{Computed $T_{2}$, $T_{5}$, $T_{7}$.  Dim~$1$.} \\ % del-longform
% \hline %% del-shortform
%  & $T_{2}$ & $(1-x)(1-2x)(1+24x+512x^{2})$ \\ % del-shortform
%  & $T_{5}$ & $(1-x)(1-5x)(1+2242x)(1-2392x)$ \\ % del-shortform
%  & $T_{7}$ & $(1-x)(1-7x)(1-1657x)(1+3617x)$ \\ % del-shortform
\hline
1 & $1$ & $\aaa \varepsilon^{0} \oplus \aaa \varepsilon^{1} \oplus \aaa \varepsilon^{2} \sigma_{3.6.\mathrm{a}.\mathrm{a}}$ \\
%% \hline
%% \multicolumn{3}{|l|}{$\sigma_{3} = \mathrm{Cusp\ form\ }q - 6q^{2} + 9q^{3} + 4q^{4} + 6q^{5} - 54q^{6} - 40q^{7} + O(q^{8})$} \\
%% \multicolumn{3}{|l|}{$\mathrm{for\ }\Gamma_0(3)\mathrm{\ with\ character\ }\aaa\mathrm{\ in\ }\mathbb{F}$} \\
%% \multicolumn{3}{|l|}{$\mathrm{of\ weight\ }6\mathrm{\ with\ coeffs\ in\ }\mathbb{Q}$} \\
\hline
\end{tabular}
\end{center}

\begin{center}
\begin{tabular}{|r|r|c|}
\hline
\multicolumn{3}{|l|}{\textbf{Level $N = 3$.  Coeffs $\Sym^5(V)\otimes\chi_3$.  Field $\mathbb{F} = GF(12379)$.}} \\
\multicolumn{3}{|l|}{Computed $T_{2}$, $T_{5}$, $T_{7,1}$.  Dim~$3$.} \\ % del-longform
% \hline %% del-shortform
%  & $T_{2}$ & $(1-x)(1-2x)(1+4x)(1-256x)$ \\ % del-shortform
%  & $T_{5}$ & $(1-x)(1-5x)(1+25x)(1+5503x)$ \\ % del-shortform
%  & $T_{7}$ & $1+3756x + O(x^2)$ \\ % del-shortform
\hline
1 & $1$ & $\aaa \varepsilon^{0} \oplus \aaa \varepsilon^{1} \oplus \chi_{3} \varepsilon^{2} \oplus \aaa \varepsilon^{8}$ \\
% \hline %% del-shortform
% \hline %% del-shortform
%  & $T_{2}$ & $(1+x)(1-2x)(1-4x)(1-256x)$ \\ % del-shortform
%  & $T_{5}$ & $(1+x)(1-5x)(1-25x)(1+5503x)$ \\ % del-shortform
%  & $T_{7}$ & $1+3756x + O(x^2)$ \\ % del-shortform
\hline
1 & $1$ & $\chi_{3} \varepsilon^{0} \oplus \aaa \varepsilon^{1} \oplus \aaa \varepsilon^{2} \oplus \aaa \varepsilon^{8}$ \\
% \hline %% del-shortform
% \hline %% del-shortform
%  & $T_{2}$ & $(1-x)(1-2x)(1+32x)(1-32x)$ \\ % del-shortform
%  & $T_{5}$ & $(1-x)(1-5x)(1+3125x)(1-3125x)$ \\ % del-shortform
%  & $T_{7}$ & $1+1627x + O(x^2)$ \\ % del-shortform
\hline
1 & $1$ & $\aaa \varepsilon^{0} \oplus \aaa \varepsilon^{1} \oplus \aaa \varepsilon^{2} \sigma_{3.7.\mathrm{b}.\mathrm{a}}$ \\
%% \hline
%% \multicolumn{3}{|l|}{$\sigma_{3} = \mathrm{Cusp\ form\ }q - 27q^{3} + 64q^{4} - 286q^{7} + O(q^{8})$} \\
%% \multicolumn{3}{|l|}{$\mathrm{for\ }\Gamma_1(3)\mathrm{\ with\ character\ }\chi_{3}\mathrm{\ in\ }\mathbb{F}$} \\
%% \multicolumn{3}{|l|}{$\mathrm{of\ weight\ }7\mathrm{\ with\ coeffs\ in\ }\mathbb{Q}$} \\
\hline
\end{tabular}
\end{center}

\begin{center}
\begin{tabular}{|r|r|c|}
\hline
\multicolumn{3}{|l|}{\textbf{Level $N = 3$.  Coeffs $\Sym^6(V)$.  Field $\mathbb{F} = GF(12379)$.}} \\
\multicolumn{3}{|l|}{Computed $T_{2}$, $T_{5}$, $T_{7,1}$.  Dim~$1$.} \\ % del-longform
% \hline %% del-shortform
%  & $T_{2}$ & $(1-x)(1-2x)(1+3923x)(1-3947x)$ \\ % del-shortform
%  & $T_{5}$ & $(1-x)(1-5x)(1+960x)(1+1669x)$ \\ % del-shortform
%  & $T_{7}$ & $1+3128x + O(x^2)$ \\ % del-shortform
\hline
1 & $1$ & $\aaa \varepsilon^{0} \oplus \aaa \varepsilon^{1} \oplus \aaa \varepsilon^{2} \sigma_{3.8.\mathrm{a}.\mathrm{a}}$ \\
%% \hline
%% \multicolumn{3}{|l|}{$\sigma_{3} = \mathrm{Cusp\ form\ }q + 6q^{2} - 27q^{3} - 92q^{4} + 390q^{5} - 162q^{6} - 64q^{7} + O(q^{8})$} \\
%% \multicolumn{3}{|l|}{$\mathrm{for\ }\Gamma_0(3)\mathrm{\ with\ character\ }\aaa\mathrm{\ in\ }\mathbb{F}$} \\
%% \multicolumn{3}{|l|}{$\mathrm{of\ weight\ }8\mathrm{\ with\ coeffs\ in\ }\mathbb{Q}$} \\
\hline
\end{tabular}
\end{center}

\begin{center}
\begin{tabular}{|r|r|c|}
\hline
\multicolumn{3}{|l|}{\textbf{Level $N = 3$.  Coeffs $\Sym^7(V)\otimes\chi_3$.  Field $\mathbb{F} = GF(12379^{2})$.}} \\
\multicolumn{3}{|l|}{Computed $T_{2}$, $T_{5,1}$, $T_{7,1}$.  Dim~$4$.} \\ % del-longform
% \hline %% del-shortform
%  & $T_{2}$ & $(1-x)(1-2x)(1+4x)(1-1024x)$ \\ % del-shortform
%  & $T_{5}$ & $1+1425x + O(x^2)$ \\ % del-shortform
%  & $T_{7}$ & $1+1095x + O(x^2)$ \\ % del-shortform
\hline
1 & $1$ & $\aaa \varepsilon^{0} \oplus \aaa \varepsilon^{1} \oplus \chi_{3} \varepsilon^{2} \oplus \aaa \varepsilon^{10}$ \\
% \hline %% del-shortform
% \hline %% del-shortform
%  & $T_{2}$ & $(1+x)(1-2x)(1-4x)(1-1024x)$ \\ % del-shortform
%  & $T_{5}$ & $1+1377x + O(x^2)$ \\ % del-shortform
%  & $T_{7}$ & $1+1095x + O(x^2)$ \\ % del-shortform
\hline
1 & $1$ & $\chi_{3} \varepsilon^{0} \oplus \aaa \varepsilon^{1} \oplus \aaa \varepsilon^{2} \oplus \aaa \varepsilon^{10}$ \\
% \hline %% del-shortform
% \hline %% del-shortform
%  & $T_{2}$ & $(1-x)(1-2x)(1+3208x-4096x^{2})$ \\ % del-shortform
%  & $T_{5}$ & $1-2442x + O(x^2)$ \\ % del-shortform
%  & $T_{7}$ & $1-911x + O(x^2)$ \\ % del-shortform
\hline
2 & $1$ & $\aaa \varepsilon^{0} \oplus \aaa \varepsilon^{1} \oplus \aaa \varepsilon^{2} \sigma_{3.9.\mathrm{b}.\mathrm{a}}$ \\
\hline
\end{tabular}
\end{center}

\subsubsection*{Level 4}

\begin{center}
\begin{tabular}{|r|r|c|}
\hline
\multicolumn{3}{|l|}{\textbf{Level $N = 4$.  Coeffs $V\otimes\chi_4$.  Field $\mathbb{F} = GF(12379)$.}} \\
\multicolumn{3}{|l|}{Computed $T_{3}$, $T_{5}$, $T_{7}$.  Dim~$2$.} \\ % del-longform
% \hline %% del-shortform
%  & $T_{3}$ & $(1-x)(1-3x)(1+9x)(1-81x)$ \\ % del-shortform
%  & $T_{5}$ & $(1-x)(1-5x)(1-25x)(1-625x)$ \\ % del-shortform
%  & $T_{7}$ & $(1-x)(1-7x)(1+49x)(1-2401x)$ \\ % del-shortform
\hline
1 & $1$ & $\aaa \varepsilon^{0} \oplus \aaa \varepsilon^{1} \oplus \chi_{4} \varepsilon^{2} \oplus \aaa \varepsilon^{4}$ \\
% \hline %% del-shortform
% \hline %% del-shortform
%  & $T_{3}$ & $(1+x)(1-3x)(1-9x)(1-81x)$ \\ % del-shortform
%  & $T_{5}$ & $(1-x)(1-5x)(1-25x)(1-625x)$ \\ % del-shortform
%  & $T_{7}$ & $(1+x)(1-7x)(1-49x)(1-2401x)$ \\ % del-shortform
\hline
1 & $1$ & $\chi_{4} \varepsilon^{0} \oplus \aaa \varepsilon^{1} \oplus \aaa \varepsilon^{2} \oplus \aaa \varepsilon^{4}$ \\
\hline
\end{tabular}
\end{center}

\begin{center}
  \begin{tabular}{|l|}
  \hline
  \textbf{Level $N = 4$.  Coeffs $\Sym^2(V)$.  Field $\mathbb{F} = GF(12379)$.} \\
  \hline
  Dim~$0$. \\
  \hline
  \end{tabular}
\end{center}

\subsubsection*{Level 5}

\begin{center}
\begin{tabular}{|r|r|c|}
\hline
\multicolumn{3}{|l|}{\textbf{Level $N = 5$.  Coeffs $V\otimes\chi_5$.  Field $\mathbb{F} = GF(16001)$.}} \\
\multicolumn{3}{|l|}{Computed $T_{2}$, $T_{3}$, $T_{7}$.  Dim~$2$.} \\ % del-longform
% \hline %% del-shortform
%  & $T_{2}$ & $(1-x)(1-2x)(1-16x)(1+2580x)$ \\ % del-shortform
%  & $T_{3}$ & $(1-x)(1-3x)(1-81x)(1-5805x)$ \\ % del-shortform
%  & $T_{7}$ & $(1-x)(1-7x)(1-397x)(1-2401x)$ \\ % del-shortform
\hline
1 & $1$ & $\aaa \varepsilon^{0} \oplus \aaa \varepsilon^{1} \oplus \chi_{5} \varepsilon^{2} \oplus \aaa \varepsilon^{4}$ \\
% \hline %% del-shortform
% \hline %% del-shortform
%  & $T_{2}$ & $(1-2x)(1-4x)(1-16x)(1+645x)$ \\ % del-shortform
%  & $T_{3}$ & $(1-3x)(1-9x)(1-81x)(1-645x)$ \\ % del-shortform
%  & $T_{7}$ & $(1-7x)(1-49x)(1+645x)(1-2401x)$ \\ % del-shortform
\hline
1 & $1$ & $\chi_{5} \varepsilon^{0} \oplus \aaa \varepsilon^{1} \oplus \aaa \varepsilon^{2} \oplus \aaa \varepsilon^{4}$ \\
\hline
\end{tabular}
\end{center}

\begin{center}
\begin{tabular}{|r|r|c|}
\hline
\multicolumn{3}{|l|}{\textbf{Level $N = 5$.  Coeffs $\Sym^2(V)$.  Field $\mathbb{F} = GF(12379^{2})$.}} \\
\multicolumn{3}{|l|}{Computed $T_{2}$, $T_{3}$, $T_{7}$, $T_{11}$.  Dim~$1$.} \\ % del-longform
% \hline %% del-shortform
%  & $T_{2}$ & $(1-x)(1-2x)(1+16x+128x^{2})$ \\ % del-shortform
%  & $T_{3}$ & $(1-x)(1-3x)(1+4297x)(1-4315x)$ \\ % del-shortform
%  & $T_{7}$ & $(1-x)(1-7x)(1-294x-5850x^{2})$ \\ % del-shortform
%  & $T_{11}$ & $(1-x)(1-11x)(1+93x)(1-3965x)$ \\ % del-shortform
\hline
1 & $1$ & $\aaa \varepsilon^{0} \oplus \aaa \varepsilon^{1} \oplus \aaa \varepsilon^{2} \sigma_{5.4.\mathrm{a}.\mathrm{a}}$ \\
%% \hline
%% \multicolumn{3}{|l|}{$\sigma_{5.4.\mathrm{a}.\mathrm{a}} = \mathrm{Cusp\ form\ }q - 4q^{2} + 2q^{3} + 8q^{4} - 5q^{5} - 8q^{6} + 6q^{7} - 23q^{9} + 20q^{10} + 32q^{11}$} \\
%% \multicolumn{3}{|c|}{$ + O(q^{12})$} \\
%% \multicolumn{3}{|l|}{$\mathrm{for\ }\Gamma_0(5)\mathrm{\ with\ character\ }\aaa\mathrm{\ in\ }\mathbb{F}$} \\
%% \multicolumn{3}{|l|}{$\mathrm{of\ weight\ }4\mathrm{\ with\ coeffs\ in\ }\mathbb{Q}$} \\
\hline
\end{tabular}
\end{center}

\begin{center}
  \begin{tabular}{|l|}
  \hline
  \textbf{Level $N = 5$.  Coeffs $\Sym^2(V)\otimes\chi_5^2$.  Field $\mathbb{F} = GF(16001)$.} \\
  \hline
  Dim~$0$. \\
  \hline
  \end{tabular}
\end{center}

\begin{center}
\begin{tabular}{|r|r|c|}
\hline
\multicolumn{3}{|l|}{\textbf{Level $N = 5$.  Coeffs $\Sym^3(V)\otimes\chi_5$.  Field $\mathbb{F} = GF(16001^{2})$.}} \\
\multicolumn{3}{|l|}{Computed $T_{2}$, $T_{3}$, $T_{7}$.  Dim~$3$.} \\ % del-longform
% \hline %% del-shortform
%  & $T_{2}$ & $(1-x)(1-2x)(1-64x)(1+2580x)$ \\ % del-shortform
%  & $T_{3}$ & $(1-x)(1-3x)(1-729x)(1-5805x)$ \\ % del-shortform
%  & $T_{7}$ & $(1-x)(1-7x)(1-397x)(1-5642x)$ \\ % del-shortform
\hline
1 & $1$ & $\aaa \varepsilon^{0} \oplus \aaa \varepsilon^{1} \oplus \chi_{5} \varepsilon^{2} \oplus \aaa \varepsilon^{6}$ \\
% \hline %% del-shortform
% \hline %% del-shortform
%  & $T_{2}$ & $(1-2x)(1-4x)(1-64x)(1+645x)$ \\ % del-shortform
%  & $T_{3}$ & $(1-3x)(1-9x)(1-645x)(1-729x)$ \\ % del-shortform
%  & $T_{7}$ & $(1-7x)(1-49x)(1+645x)(1-5642x)$ \\ % del-shortform
\hline
1 & $1$ & $\chi_{5} \varepsilon^{0} \oplus \aaa \varepsilon^{1} \oplus \aaa \varepsilon^{2} \oplus \aaa \varepsilon^{6}$ \\
% \hline %% del-shortform
% \hline %% del-shortform
%  & $T_{2}$ & $(1-x)(1-2x)(1+3906x)(1-6482x)$ \\ % del-shortform
%  & $T_{3}$ & $(1-x)(1-3x)(1+2882x+7581x^{2})$ \\ % del-shortform
%  & $T_{7}$ & $(1-x)(1-7x)(1-4405x-266x^{2})$ \\ % del-shortform
\hline
1 & $1$ & $\aaa \varepsilon^{0} \oplus \aaa \varepsilon^{1} \oplus \aaa \varepsilon^{2} \sigma_{5.5.\mathrm{c}.\mathrm{a}}$ \\
%% \hline
%% \multicolumn{3}{|l|}{$\sigma_{5} = \mathrm{Cusp\ form\ }q + \left(-\zeta_{4} - 1\right)q^{2} + \left(6 \zeta_{4} - 6\right)q^{3} - 14 \zeta_{4}q^{4}$} \\
%% \multicolumn{3}{|c|}{$ + \left(15 \zeta_{4} + 20\right)q^{5} + 12q^{6} + \left(-26 \zeta_{4} - 26\right)q^{7} + O(q^{8})$} \\
%% \multicolumn{3}{|l|}{$\mathrm{for\ }\Gamma_1(5)\mathrm{\ with\ character\ }\chi_{5}\mathrm{\ in\ }\mathbb{F}$} \\
%% \multicolumn{3}{|l|}{$\mathrm{of\ weight\ }5\mathrm{\ with\ coeffs\ in\ }\mathbb{Q}(\zeta_{4})$} \\
\hline
\end{tabular}
\end{center}

\begin{center}
\begin{tabular}{|r|r|c|}
\hline
\multicolumn{3}{|l|}{\textbf{Level $N = 5$.  Coeffs $\Sym^4(V)$.  Field $\mathbb{F} = GF(16001^{2})$.}} \\
\multicolumn{3}{|l|}{Computed $T_{2}$, $T_{3}$, $T_{7}$.  Dim~$1$.} \\ % del-longform
% \hline %% del-shortform
%  & $T_{2}$ & $(1-x)(1-2x)(1+4882x)(1-4890x)$ \\ % del-shortform
%  & $T_{3}$ & $(1-x)(1-3x)(1+36x+3682x^{2})$ \\ % del-shortform
%  & $T_{7}$ & $(1-x)(1-7x)(1-1616x)(1-7792x)$ \\ % del-shortform
\hline
1 & $1$ & $\aaa \varepsilon^{0} \oplus \aaa \varepsilon^{1} \oplus \aaa \varepsilon^{2} \sigma_{5.6.\mathrm{a}.\mathrm{a}}$ \\
%% \hline
%% \multicolumn{3}{|l|}{$\sigma_{5} = \mathrm{Cusp\ form\ }q + 2q^{2} - 4q^{3} - 28q^{4} + 25q^{5} - 8q^{6} + 192q^{7} + O(q^{8})$} \\
%% \multicolumn{3}{|l|}{$\mathrm{for\ }\Gamma_0(5)\mathrm{\ with\ character\ }\aaa\mathrm{\ in\ }\mathbb{F}$} \\
%% \multicolumn{3}{|l|}{$\mathrm{of\ weight\ }6\mathrm{\ with\ coeffs\ in\ }\mathbb{Q}$} \\
\hline
\end{tabular}
\end{center}

\begin{center}
\begin{tabular}{|r|r|c|}
\hline
\multicolumn{3}{|l|}{\textbf{Level $N = 5$.  Coeffs $\Sym^4(V)\otimes\chi_5^2$.  Field $\mathbb{F} = GF(16001^{2})$.}} \\
\multicolumn{3}{|l|}{Computed $T_{2}$, $T_{3}$, $T_{7}$.  Dim~$2$.} \\ % del-longform
% \hline %% del-shortform
%  & $T_{2}$ & $(1-x)(1-2x)(1 + \left(15993 a + 570\right)x)(1 + \left(15993 a + 15439\right)x) \mathrm{\ with\ } a^2 + 16000*a + 3 = 0$ \\ % del-shortform
%  & $T_{3}$ & $(1-x)(1-3x)(1 + \left(54 a + 81\right)x)(1 + \left(54 a + 15866\right)x) \mathrm{\ with\ } a^2 + 16000*a + 3 = 0$ \\ % del-shortform
%  & $T_{7}$ & $(1-x)(1-7x)(1 + \left(882 a + 2077\right)x)(1 + \left(882 a + 13042\right)x) \mathrm{\ with\ } a^2 + 16000*a + 3 = 0$ \\ % del-shortform
\hline
2 & $1$ & $\aaa \varepsilon^{0} \oplus \aaa \varepsilon^{1} \oplus \aaa \varepsilon^{2} \sigma_{5.6.\mathrm{b}.\mathrm{a}}$ \\
\hline
\end{tabular}
\end{center}

\begin{center}
\begin{tabular}{|r|r|c|}
\hline
\multicolumn{3}{|l|}{\textbf{Level $N = 5$.  Coeffs $\Sym^5(V)\otimes\chi_5$.  Field $\mathbb{F} = GF(16001^{6})$.}} \\
\multicolumn{3}{|l|}{Computed $T_{2}$, $T_{3}$, $T_{7,1}$.  Dim~$4$.} \\ % del-longform
% \hline %% del-shortform
%  & $T_{2}$ & $(1-x)(1-2x)(1-256x)(1+2580x)$ \\ % del-shortform
%  & $T_{3}$ & $(1-x)(1-3x)(1-5805x)(1-6561x)$ \\ % del-shortform
%  & $T_{7}$ & $1-4846x + O(x^2)$ \\ % del-shortform
\hline
1 & $1$ & $\aaa \varepsilon^{0} \oplus \aaa \varepsilon^{1} \oplus \chi_{5}^{3} \varepsilon^{2} \oplus \aaa \varepsilon^{8}$ \\
% \hline %% del-shortform
% \hline %% del-shortform
%  & $T_{2}$ & $(1-2x)(1-4x)(1-256x)(1+645x)$ \\ % del-shortform
%  & $T_{3}$ & $(1-3x)(1-9x)(1-645x)(1-6561x)$ \\ % del-shortform
%  & $T_{7}$ & $1-3852x + O(x^2)$ \\ % del-shortform
\hline
1 & $1$ & $\chi_{5}^{3} \varepsilon^{0} \oplus \aaa \varepsilon^{1} \oplus \aaa \varepsilon^{2} \oplus \aaa \varepsilon^{8}$ \\
% \hline %% del-shortform
% \hline %% del-shortform
%  & $T_{2}$ & $(1 + \left(13676 c + 14626\right)x)(1 + \left(2081 c + 8665\right)x + \left(3114 c + 5319\right)x^{2} + \left(7294 c + 92\right)x^{3}) \mathrm{\ with\ } c^2 + 16000*c + 3 = 0$ \\ % del-shortform
%  & $T_{3}$ & $(1 + \left(14736 c + 9472\right)x + \left(3726 c + 683\right)x^{2} + \left(4176 c + 7268\right)x^{3}+5935x^{4}) \mathrm{\ with\ } c^2 + 16000*c + 3 = 0$ \\ % del-shortform
%  & $T_{7}$ & $1 + \left(353 c + 5399\right)x + O(x^2) \mathrm{\ with\ } c^2 + 16000*c + 3 = 0$ \\ % del-shortform
\hline
2 & $1$ & $\aaa \varepsilon^{0} \oplus \aaa \varepsilon^{1} \oplus \aaa \varepsilon^{2} \sigma_{5.7.\mathrm{c}.\mathrm{a}}$ \\
\hline
\end{tabular}
\end{center}

\begin{center}
\begin{tabular}{|r|r|c|}
\hline
\multicolumn{3}{|l|}{\textbf{Level $N = 5$.  Coeffs $\Sym^6(V)$.  Field $\mathbb{F} = GF(16001^{6})$.}} \\
\multicolumn{3}{|l|}{Computed $T_{2}$, $T_{3}$, $T_{7,1}$.  Dim~$3$.} \\ % del-longform
% \hline %% del-shortform
%  & $T_{2}$ & $(1-x)(1-2x)(1+56x+2048x^{2})$ \\ % del-shortform
%  & $T_{3}$ & $(1-x)(1-3x)(1-4987x)(1+5419x)$ \\ % del-shortform
%  & $T_{7}$ & $1+543x + O(x^2)$ \\ % del-shortform
\hline
1 & $1$ & $\aaa \varepsilon^{0} \oplus \aaa \varepsilon^{1} \oplus \aaa \varepsilon^{2} \sigma_{5.8.\mathrm{a}.\mathrm{a}}$ \\
%% \hline
%% \multicolumn{3}{|l|}{$\sigma_{5} = \mathrm{Cusp\ form\ }q - 14q^{2} - 48q^{3} + 68q^{4} + 125q^{5} + 672q^{6} - 1644q^{7} + O(q^{8})$} \\
%% \multicolumn{3}{|l|}{$\mathrm{for\ }\Gamma_0(5)\mathrm{\ with\ character\ }\aaa\mathrm{\ in\ }\mathbb{F}$} \\
%% \multicolumn{3}{|l|}{$\mathrm{of\ weight\ }8\mathrm{\ with\ coeffs\ in\ }\mathbb{Q}$} \\
% \hline %% del-shortform
% \hline %% del-shortform
%  & $T_{2}$ & $(1 + \left(15067 c + 3190\right)x + \left(4459 c + 11955\right)x^{2} + \left(663 c + 6169\right)x^{3}+7195x^{4}) \mathrm{\ with\ } c^2 + 16000*c + 3 = 0$ \\ % del-shortform
%  & $T_{3}$ & $(1 + \left(1929 c + 10622\right)x)(1 + \left(14893 c + 11843\right)x + \left(14900 c + 3956\right)x^{2} + \left(7641 c + 5852\right)x^{3}) \mathrm{\ with\ } c^2 + 16000*c + 3 = 0$ \\ % del-shortform
%  & $T_{7}$ & $1 + \left(3818 c + 7658\right)x + O(x^2) \mathrm{\ with\ } c^2 + 16000*c + 3 = 0$ \\ % del-shortform
\hline
2 & $1$ & $\aaa \varepsilon^{0} \oplus \aaa \varepsilon^{1} \oplus \aaa \varepsilon^{2} \sigma_{5.8.\mathrm{a}.\mathrm{b}}$ \\
\hline
\end{tabular}
\end{center}

\begin{center}
\begin{tabular}{|r|r|c|}
\hline
\multicolumn{3}{|l|}{\textbf{Level $N = 5$.  Coeffs $\Sym^6(V)\otimes\chi_5^2$.  Field $\mathbb{F} = GF(16001^{6})$.}} \\
\multicolumn{3}{|l|}{Computed $T_{2}$, $T_{3}$, $T_{7,1}$.  Dim~$2$.} \\ % del-longform
% \hline %% del-shortform
%  & $T_{2}$ & $(1-x)(1-2x)(1+4959x-2048x^{2})$ \\ % del-shortform
%  & $T_{3}$ & $(1-x)(1-3x)(1-2529x-1136x^{2})$ \\ % del-shortform
%  & $T_{7}$ & $1+2978x + O(x^2)$ \\ % del-shortform
\hline
2 & $1$ & $\aaa \varepsilon^{0} \oplus \aaa \varepsilon^{1} \oplus \aaa \varepsilon^{2} \sigma_{5.8.\mathrm{b}.\mathrm{a}}$ \\
%% \hline
%% \multicolumn{3}{|l|}{$\sigma_{5} = \mathrm{Cusp\ form\ }q + b_{0}q^{2} + 3 b_{0}q^{3} + 12q^{4} + \left(-25 b_{0} + 75\right)q^{5} - 348q^{6} - 39 b_{0}q^{7}$} \\
%% \multicolumn{3}{|c|}{$ + O(q^{8})$} \\
%% \multicolumn{3}{|l|}{$\mathrm{for\ }\Gamma_1(5)\mathrm{\ with\ character\ }\chi_{5}^{2}\mathrm{\ in\ }\mathbb{F}$} \\
%% \multicolumn{3}{|l|}{$\mathrm{of\ weight\ }8\mathrm{\ with\ coeffs\ in\ }\mathbb{Q}[b_{0}]/(b_{0}^{2} + 116)$} \\
% \hline %% del-shortform
% \hline %% del-shortform
%  & $T_{2}$ & $(1-x)(1-2x)(1-4959x-2048x^{2})$ \\ % del-shortform
%  & $T_{3}$ & $(1-x)(1-3x)(1+2529x-1136x^{2})$ \\ % del-shortform
%  & $T_{7}$ & $1-2994x + O(x^2)$ \\ % del-shortform
%% \hline
%% 1 & $1$ & $\aaa \varepsilon^{0} \oplus \aaa \varepsilon^{1} \oplus \aaa \varepsilon^{2} \sigma_{5.8.\mathrm{b}.\mathrm{a}}$ \\
%% \hline
%% \multicolumn{3}{|l|}{$\sigma_{5} = \mathrm{Cusp\ form\ }q + b_{0}q^{2} + 3 b_{0}q^{3} + 12q^{4} + \left(-25 b_{0} + 75\right)q^{5} - 348q^{6} - 39 b_{0}q^{7}$} \\
%% \multicolumn{3}{|c|}{$ + O(q^{8})$} \\
%% \multicolumn{3}{|l|}{$\mathrm{for\ }\Gamma_1(5)\mathrm{\ with\ character\ }\chi_{5}^{2}\mathrm{\ in\ }\mathbb{F}$} \\
%% \multicolumn{3}{|l|}{$\mathrm{of\ weight\ }8\mathrm{\ with\ coeffs\ in\ }\mathbb{Q}[b_{0}]/(b_{0}^{2} + 116)$} \\
\hline
\end{tabular}
\end{center}

\begin{center}
\begin{tabular}{|r|r|c|}
\hline
\multicolumn{3}{|l|}{\textbf{Level $N = 5$.  Coeffs $\Sym^7(V)\otimes\chi_5$.  Field $\mathbb{F} = GF(16001^{6})$.}} \\
\multicolumn{3}{|l|}{Computed $T_{2}$, $T_{3}$.  Dim~$5$.} \\ % del-longform
% \hline %% del-shortform
%  & $T_{2}$ & $(1-x)(1-2x)(1-1024x)(1+2580x)$ \\ % del-shortform
%  & $T_{3}$ & $(1-x)(1-3x)(1+4955x)(1-5805x)$ \\ % del-shortform
\hline
1 & $1$ & $\aaa \varepsilon^{0} \oplus \aaa \varepsilon^{1} \oplus \chi_{5}^{3} \varepsilon^{2} \oplus \aaa \varepsilon^{10}$ \\
% \hline %% del-shortform
% \hline %% del-shortform
%  & $T_{2}$ & $(1-2x)(1-4x)(1+645x)(1-1024x)$ \\ % del-shortform
%  & $T_{3}$ & $(1-3x)(1-9x)(1-645x)(1+4955x)$ \\ % del-shortform
\hline
1 & $1$ & $\chi_{5}^{3} \varepsilon^{0} \oplus \aaa \varepsilon^{1} \oplus \aaa \varepsilon^{2} \oplus \aaa \varepsilon^{10}$ \\
% \hline %% del-shortform
% \hline %% del-shortform
%  & $T_{2}$ & $(1 + \left(14830 c^{2} + 12929 c + 6542\right)x)(1 + \left(3974 c^{2} + c + 11626\right)x + \left(12146 c^{2} + 4983 c + 5982\right)x^{2} + \left(14399 c^{2} + 749 c + 7577\right)x^{3}) \mathrm{\ with\ } c^3 + 3*c + 15998 = 0$ \\ % del-shortform
%  & $T_{3}$ & $(1 + \left(15453 c^{2} + 14852 c + 5355\right)x)(1 + \left(15505 c^{2} + 15805 c + 5547\right)x)(1 + \left(11287 c^{2} + 12265 c + 11049\right)x)(1 + \left(10636 c^{2} + 5787 c + 13461\right)x) \mathrm{\ with\ } c^3 + 3*c + 15998 = 0$ \\ % del-shortform
\hline
3 & $1$ & $\aaa \varepsilon^{0} \oplus \aaa \varepsilon^{1} \oplus \aaa \varepsilon^{2} \sigma_{5.9.\mathrm{c}.\mathrm{a}}$ \\
\hline
\end{tabular}
\end{center}

\subsubsection*{Level 6}

\begin{center}
\begin{tabular}{|r|r|c|}
\hline
\multicolumn{3}{|l|}{\textbf{Level $N = 6$.  Coeffs $V\otimes\chi_6$.  Field $\mathbb{F} = GF(12379)$.}} \\
\multicolumn{3}{|l|}{Computed $T_{5}$, $T_{7}$.  Dim~$6$.} \\ % del-longform
% \hline %% del-shortform
%  & $T_{5}$ & $(1-x)(1-5x)(1+25x)(1-625x)$ \\ % del-shortform
%  & $T_{7}$ & $(1-x)(1-7x)(1-49x)(1-2401x)$ \\ % del-shortform
\hline
1 & $3$ & $\aaa \varepsilon^{0} \oplus \aaa \varepsilon^{1} \oplus \chi_{6} \varepsilon^{2} \oplus \aaa \varepsilon^{4}$ \\
% \hline %% del-shortform
% \hline %% del-shortform
%  & $T_{5}$ & $(1+x)(1-5x)(1-25x)(1-625x)$ \\ % del-shortform
%  & $T_{7}$ & $(1-x)(1-7x)(1-49x)(1-2401x)$ \\ % del-shortform
\hline
1 & $3$ & $\chi_{6} \varepsilon^{0} \oplus \aaa \varepsilon^{1} \oplus \aaa \varepsilon^{2} \oplus \aaa \varepsilon^{4}$ \\
\hline
\end{tabular}
\end{center}

\begin{center}
\begin{tabular}{|r|r|c|}
\hline
\multicolumn{3}{|l|}{\textbf{Level $N = 6$.  Coeffs $\Sym^2(V)$.  Field $\mathbb{F} = GF(12379)$.}} \\
\multicolumn{3}{|l|}{Computed $T_{5}$, $T_{7}$, $T_{11}$.  Dim~$1$.} \\ % del-longform
% \hline %% del-shortform
%  & $T_{5}$ & $(1-x)(1-5x)(1-150x+3851x^{2})$ \\ % del-shortform
%  & $T_{7}$ & $(1-x)(1-7x)(1-4998x)(1+5782x)$ \\ % del-shortform
%  & $T_{11}$ & $(1-x)(1-11x)(1-1452x+2625x^{2})$ \\ % del-shortform
\hline
1 & $1$ & $\aaa \varepsilon^{0} \oplus \aaa \varepsilon^{1} \oplus \aaa \varepsilon^{2} \sigma_{6.4.\mathrm{a}.\mathrm{a}}$ \\
%% \hline
%% \multicolumn{3}{|l|}{$\sigma_{6} = \mathrm{Cusp\ form\ }q - 2q^{2} - 3q^{3} + 4q^{4} + 6q^{5} + 6q^{6} - 16q^{7} - 8q^{8} + 9q^{9} - 12q^{10} + 12q^{11}$} \\
%% \multicolumn{3}{|c|}{$ + O(q^{12})$} \\
%% \multicolumn{3}{|l|}{$\mathrm{for\ }\Gamma_0(6)\mathrm{\ with\ character\ }\aaa\mathrm{\ in\ }\mathbb{F}$} \\
%% \multicolumn{3}{|l|}{$\mathrm{of\ weight\ }4\mathrm{\ with\ coeffs\ in\ }\mathbb{Q}$} \\
\hline
\end{tabular}
\end{center}

\subsubsection*{Level 7}

\begin{center}
\begin{tabular}{|r|r|c|}
\hline
\multicolumn{3}{|l|}{\textbf{Level $N = 7$.  Coeffs $V\otimes\chi_7$.  Field $\mathbb{F} = GF(12037^{2})$.}} \\
\multicolumn{3}{|l|}{Computed $T_{2}$, $T_{3}$, $T_{5}$.  Dim~$2$.} \\ % del-longform
% \hline %% del-shortform
%  & $T_{2}$ & $(1-x)(1-2x)(1-16x)(1+5176x)$ \\ % del-shortform
%  & $T_{3}$ & $(1-x)(1-3x)(1-81x)(1-400x)$ \\ % del-shortform
%  & $T_{5}$ & $(1-x)(1-5x)(1-625x)(1+3761x)$ \\ % del-shortform
\hline
1 & $1$ & $\aaa \varepsilon^{0} \oplus \aaa \varepsilon^{1} \oplus \chi_{7} \varepsilon^{2} \oplus \aaa \varepsilon^{4}$ \\
% \hline %% del-shortform
% \hline %% del-shortform
%  & $T_{2}$ & $(1-2x)(1-4x)(1-16x)(1+1294x)$ \\ % del-shortform
%  & $T_{3}$ & $(1-3x)(1-9x)(1-81x)(1+1293x)$ \\ % del-shortform
%  & $T_{5}$ & $(1-5x)(1-25x)(1-625x)(1-1294x)$ \\ % del-shortform
\hline
1 & $1$ & $\chi_{7} \varepsilon^{0} \oplus \aaa \varepsilon^{1} \oplus \aaa \varepsilon^{2} \oplus \aaa \varepsilon^{4}$ \\
\hline
\end{tabular}
\end{center}

\begin{center}
\begin{tabular}{|r|r|c|}
\hline
\multicolumn{3}{|l|}{\textbf{Level $N = 7$.  Coeffs $V\otimes\chi_7^3$.  Field $\mathbb{F} = GF(12037^{2})$.}} \\
\multicolumn{3}{|l|}{Computed $T_{2}$, $T_{3}$, $T_{5}$.  Dim~$4$.} \\ % del-longform
% \hline %% del-shortform
%  & $T_{2}$ & $(1-x)(1-2x)(1-4x)(1-16x)$ \\ % del-shortform
%  & $T_{3}$ & $(1-x)(1-3x)(1+9x)(1-81x)$ \\ % del-shortform
%  & $T_{5}$ & $(1-x)(1-5x)(1+25x)(1-625x)$ \\ % del-shortform
\hline
1 & $1$ & $\aaa \varepsilon^{0} \oplus \aaa \varepsilon^{1} \oplus \chi_{7}^{3} \varepsilon^{2} \oplus \aaa \varepsilon^{4}$ \\
% \hline %% del-shortform
% \hline %% del-shortform
%  & $T_{2}$ & $(1-x)(1-2x)(1-4x)(1-16x)$ \\ % del-shortform
%  & $T_{3}$ & $(1+x)(1-3x)(1-9x)(1-81x)$ \\ % del-shortform
%  & $T_{5}$ & $(1+x)(1-5x)(1-25x)(1-625x)$ \\ % del-shortform
\hline
1 & $1$ & $\chi_{7}^{3} \varepsilon^{0} \oplus \aaa \varepsilon^{1} \oplus \aaa \varepsilon^{2} \oplus \aaa \varepsilon^{4}$ \\
% \hline %% del-shortform
% \hline %% del-shortform
%  & $T_{2}$ & $(1-x)(1-2x)(1-2191x)(1+2203x)$ \\ % del-shortform
%  & $T_{3}$ & $(1-x)(1-3x)(1+27x)(1-27x)$ \\ % del-shortform
%  & $T_{5}$ & $(1-x)(1-5x)(1+125x)(1-125x)$ \\ % del-shortform
\hline
1 & $1$ & $\aaa \varepsilon^{0} \oplus \aaa \varepsilon^{1} \oplus \aaa \varepsilon^{2} \sigma_{7.3.\mathrm{b}.\mathrm{a}}$ \\
\hline
\end{tabular}
\end{center}

\begin{center}
\begin{tabular}{|r|r|c|}
\hline
\multicolumn{3}{|l|}{\textbf{Level $N = 7$.  Coeffs $\Sym^2(V)$.  Field $\mathbb{F} = GF(12379^{2})$.}} \\
\multicolumn{3}{|l|}{Computed $T_{2}$, $T_{3}$, $T_{5}$, $T_{11}$.  Dim~$1$.} \\ % del-longform
% \hline %% del-shortform
%  & $T_{2}$ & $(1-x)(1-2x)(1-4051x)(1+4055x)$ \\ % del-shortform
%  & $T_{3}$ & $(1-x)(1-3x)(1-4297x)(1+4315x)$ \\ % del-shortform
%  & $T_{5}$ & $(1-x)(1-5x)(1-400x+3851x^{2})$ \\ % del-shortform
%  & $T_{11}$ & $(1-x)(1-11x)(1-162x)(1+1130x)$ \\ % del-shortform
\hline
1 & $1$ & $\aaa \varepsilon^{0} \oplus \aaa \varepsilon^{1} \oplus \aaa \varepsilon^{2} \sigma_{7.4.\mathrm{a}.\mathrm{a}}$ \\
%% \hline
%% \multicolumn{3}{|l|}{$\sigma_{7} = \mathrm{Cusp\ form\ }q - q^{2} - 2q^{3} - 7q^{4} + 16q^{5} + 2q^{6} - 7q^{7} + 15q^{8} - 23q^{9} - 16q^{10}$} \\
%% \multicolumn{3}{|c|}{$ - 8q^{11} + O(q^{12})$} \\
%% \multicolumn{3}{|l|}{$\mathrm{for\ }\Gamma_0(7)\mathrm{\ with\ character\ }\aaa\mathrm{\ in\ }\mathbb{F}$} \\
%% \multicolumn{3}{|l|}{$\mathrm{of\ weight\ }4\mathrm{\ with\ coeffs\ in\ }\mathbb{Q}$} \\
\hline
\end{tabular}
\end{center}

\begin{center}
\begin{tabular}{|r|r|c|}
\hline
\multicolumn{3}{|l|}{\textbf{Level $N = 7$.  Coeffs $\Sym^2(V)\otimes\chi_7^2$.  Field $\mathbb{F} = GF(12037^{2})$.}} \\
\multicolumn{3}{|l|}{Computed $T_{2}$, $T_{3}$, $T_{5}$.  Dim~$1$.} \\ % del-longform
% \hline %% del-shortform
%  & $T_{2}$ & $(1-x)(1-2x)(1+804x)(1-2489x)$ \\ % del-shortform
%  & $T_{3}$ & $(1-x)(1-3x)(1+962x)(1+1838x)$ \\ % del-shortform
%  & $T_{5}$ & $(1-x)(1-5x)(1-2253x+1121x^{2})$ \\ % del-shortform
\hline
1 & $1$ & $\aaa \varepsilon^{0} \oplus \aaa \varepsilon^{1} \oplus \aaa \varepsilon^{2} \sigma_{7.4.\mathrm{c}.\mathrm{a}}$ \\
%% \hline
%% \multicolumn{3}{|l|}{$\sigma_{7} = \mathrm{Cusp\ form\ }q + \left(2 \zeta_{6} - 2\right)q^{2} - 7 \zeta_{6}q^{3} + 4 \zeta_{6}q^{4} + \left(7 \zeta_{6} - 7\right)q^{5}$} \\
%% \multicolumn{3}{|c|}{$ + 14q^{6} + \left(-14 \zeta_{6} + 21\right)q^{7} + O(q^{8})$} \\
%% \multicolumn{3}{|l|}{$\mathrm{for\ }\Gamma_1(7)\mathrm{\ with\ character\ }\chi_{7}^{2}\mathrm{\ in\ }\mathbb{F}$} \\
%% \multicolumn{3}{|l|}{$\mathrm{of\ weight\ }4\mathrm{\ with\ coeffs\ in\ }\mathbb{Q}(\zeta_{6})$} \\
\hline
\end{tabular}
\end{center}

\begin{center}
\begin{tabular}{|r|r|c|}
\hline
\multicolumn{3}{|l|}{\textbf{Level $N = 7$.  Coeffs $\Sym^3(V)\otimes\chi_7$.  Field $\mathbb{F} = GF(12037^{2})$.}} \\
\multicolumn{3}{|l|}{Computed $T_{2}$, $T_{3}$, $T_{5}$.  Dim~$4$.} \\ % del-longform
% \hline %% del-shortform
%  & $T_{2}$ & $(1-x)(1-2x)(1-64x)(1+5176x)$ \\ % del-shortform
%  & $T_{3}$ & $(1-x)(1-3x)(1-400x)(1-729x)$ \\ % del-shortform
%  & $T_{5}$ & $(1-x)(1-5x)(1-3588x)(1+3761x)$ \\ % del-shortform
\hline
1 & $1$ & $\aaa \varepsilon^{0} \oplus \aaa \varepsilon^{1} \oplus \chi_{7} \varepsilon^{2} \oplus \aaa \varepsilon^{6}$ \\
% \hline %% del-shortform
% \hline %% del-shortform
%  & $T_{2}$ & $(1-2x)(1-4x)(1-64x)(1+1294x)$ \\ % del-shortform
%  & $T_{3}$ & $(1-3x)(1-9x)(1-729x)(1+1293x)$ \\ % del-shortform
%  & $T_{5}$ & $(1-5x)(1-25x)(1-1294x)(1-3588x)$ \\ % del-shortform
\hline
1 & $1$ & $\chi_{7} \varepsilon^{0} \oplus \aaa \varepsilon^{1} \oplus \aaa \varepsilon^{2} \oplus \aaa \varepsilon^{6}$ \\
% \hline %% del-shortform
% \hline %% del-shortform
%  & $T_{2}$ & $(1-x)(1-2x)(1 + \left(11053 a + 6947\right)x)(1 + \left(4363 a + 11112\right)x) \mathrm{\ with\ } a^2 + 12036*a + 5 = 0$ \\ % del-shortform
%  & $T_{3}$ & $(1-x)(1-3x)(1 + \left(8936 a + 4393\right)x)(1 + \left(6815 a + 5378\right)x) \mathrm{\ with\ } a^2 + 12036*a + 5 = 0$ \\ % del-shortform
%  & $T_{5}$ & $(1-x)(1-5x)(1 + \left(7842 a + 321\right)x)(1 + \left(1650 a + 327\right)x) \mathrm{\ with\ } a^2 + 12036*a + 5 = 0$ \\ % del-shortform
\hline
2 & $1$ & $\aaa \varepsilon^{0} \oplus \aaa \varepsilon^{1} \oplus \aaa \varepsilon^{2} \sigma_{7.5.\mathrm{d}.\mathrm{a}}$ \\
\hline
\end{tabular}
\end{center}

\begin{center}
\begin{tabular}{|r|r|c|}
\hline
\multicolumn{3}{|l|}{\textbf{Level $N = 7$.  Coeffs $\Sym^3(V)\otimes\chi_7^3$.  Field $\mathbb{F} = GF(12037^{2})$.}} \\
\multicolumn{3}{|l|}{Computed $T_{2}$, $T_{3}$, $T_{5}$.  Dim~$4$.} \\ % del-longform
% \hline %% del-shortform
%  & $T_{2}$ & $(1-x)(1-2x)(1-4x)(1-64x)$ \\ % del-shortform
%  & $T_{3}$ & $(1-x)(1-3x)(1+9x)(1-729x)$ \\ % del-shortform
%  & $T_{5}$ & $(1-x)(1-5x)(1+25x)(1-3588x)$ \\ % del-shortform
\hline
1 & $1$ & $\aaa \varepsilon^{0} \oplus \aaa \varepsilon^{1} \oplus \chi_{7}^{3} \varepsilon^{2} \oplus \aaa \varepsilon^{6}$ \\
% \hline %% del-shortform
% \hline %% del-shortform
%  & $T_{2}$ & $(1-x)(1-2x)(1-4x)(1-64x)$ \\ % del-shortform
%  & $T_{3}$ & $(1+x)(1-3x)(1-9x)(1-729x)$ \\ % del-shortform
%  & $T_{5}$ & $(1+x)(1-5x)(1-25x)(1-3588x)$ \\ % del-shortform
\hline
1 & $1$ & $\chi_{7}^{3} \varepsilon^{0} \oplus \aaa \varepsilon^{1} \oplus \aaa \varepsilon^{2} \oplus \aaa \varepsilon^{6}$ \\
% \hline %% del-shortform
% \hline %% del-shortform
%  & $T_{2}$ & $(1-2x)(1-64x)(1-3557x)(1+3560x)$ \\ % del-shortform
%  & $T_{3}$ & $(1+3x)(1-3x)^2(1-729x)$ \\ % del-shortform
%  & $T_{5}$ & $(1+5x)(1-5x)^2(1-3588x)$ \\ % del-shortform
\hline
1 & $1$ & $\aaa \varepsilon^{1} \oplus \aaa \varepsilon^{6} \oplus \aaa \varepsilon^{0} \sigma_{7.3.\mathrm{b}.\mathrm{a}}$ \\
\hline
\end{tabular}
\end{center}

\begin{center}
\begin{tabular}{|r|r|c|}
\hline
\multicolumn{3}{|l|}{\textbf{Level $N = 7$.  Coeffs $\Sym^4(V)$.  Field $\mathbb{F} = GF(12037^{6})$.}} \\
\multicolumn{3}{|l|}{Computed $T_{2}$, $T_{3}$, $T_{5}$.  Dim~$3$.} \\ % del-longform
% \hline %% del-shortform
%  & $T_{2}$ & $(1-x)(1-2x)(1-4374x)(1+4414x)$ \\ % del-shortform
%  & $T_{3}$ & $(1-x)(1-3x)(1+126x-4391x^{2})$ \\ % del-shortform
%  & $T_{5}$ & $(1-x)(1-5x)(1-153x)(1+1553x)$ \\ % del-shortform
\hline
1 & $1$ & $\aaa \varepsilon^{0} \oplus \aaa \varepsilon^{1} \oplus \aaa \varepsilon^{2} \sigma_{7.6.\mathrm{a}.\mathrm{a}}$ \\
%% \hline
%% \multicolumn{3}{|l|}{$\sigma_{7} = \mathrm{Cusp\ form\ }q - 10q^{2} - 14q^{3} + 68q^{4} - 56q^{5} + 140q^{6} - 49q^{7} + O(q^{8})$} \\
%% \multicolumn{3}{|l|}{$\mathrm{for\ }\Gamma_0(7)\mathrm{\ with\ character\ }\aaa\mathrm{\ in\ }\mathbb{F}$} \\
%% \multicolumn{3}{|l|}{$\mathrm{of\ weight\ }6\mathrm{\ with\ coeffs\ in\ }\mathbb{Q}$} \\
% \hline %% del-shortform
% \hline %% del-shortform
%  & $T_{2}$ & $(1 + \left(3811 c + 8717\right)x)(1 + \left(7148 c + 7972\right)x + \left(3814 c + 916\right)x^{2} + \left(9179 c + 10687\right)x^{3}) \mathrm{\ with\ } c^2 + 12036*c + 5 = 0$ \\ % del-shortform
%  & $T_{3}$ & $(1 + \left(6338 c + 1040\right)x)(1 + \left(6168 c + 9318\right)x + \left(5504 c + 4040\right)x^{2} + \left(5949 c + 8424\right)x^{3}) \mathrm{\ with\ } c^2 + 12036*c + 5 = 0$ \\ % del-shortform
%  & $T_{5}$ & $(1 + \left(2213 c + 4415\right)x + \left(10442 c + 5238\right)x^{2})(1 + \left(7773 c + 347\right)x + \left(10672 c + 10806\right)x^{2}) \mathrm{\ with\ } c^2 + 12036*c + 5 = 0$ \\ % del-shortform
\hline
2 & $1$ & $\aaa \varepsilon^{0} \oplus \aaa \varepsilon^{1} \oplus \aaa \varepsilon^{2} \sigma_{7.6.\mathrm{a}.\mathrm{b}}$ \\
\hline
\end{tabular}
\end{center}

\begin{center}
\begin{tabular}{|r|r|c|}
\hline
\multicolumn{3}{|l|}{\textbf{Level $N = 7$.  Coeffs $\Sym^4(V)\otimes\chi_7^2$.  Field $\mathbb{F} = GF(12037^{6})$.}} \\
\multicolumn{3}{|l|}{Computed $T_{2}$, $T_{3}$, $T_{5}$.  Dim~$2$.} \\ % del-longform
% \hline %% del-shortform
%  & $T_{2}$ & $(1-x)(1-2x)(1+1396x)(1+4682x)$ \\ % del-shortform
%  & $T_{3}$ & $(1-x)(1-3x)(1+507x+490x^{2})$ \\ % del-shortform
%  & $T_{5}$ & $(1-x)(1-5x)(1+3047x+3951x^{2})$ \\ % del-shortform
\hline
2 & $1$ & $\aaa \varepsilon^{0} \oplus \aaa \varepsilon^{1} \oplus \aaa \varepsilon^{2} \sigma_{7.6.\mathrm{c}.\mathrm{a}}$ \\
\hline
\end{tabular}
\end{center}

\begin{center}
\begin{tabular}{|r|r|c|}
\hline
\multicolumn{3}{|l|}{\textbf{Level $N = 7$.  Coeffs $\Sym^5(V)\otimes\chi_7$.  Field $\mathbb{F} = GF(12037^{6})$.}} \\
\multicolumn{3}{|l|}{Computed $T_{2}$, $T_{3}$, $T_{5,1}$.  Dim~$5$.} \\ % del-longform
% \hline %% del-shortform
%  & $T_{2}$ & $(1-x)(1-2x)(1-256x)(1+5176x)$ \\ % del-shortform
%  & $T_{3}$ & $(1-x)(1-3x)(1-400x)(1+5476x)$ \\ % del-shortform
%  & $T_{5}$ & $1-1686x + O(x^2)$ \\ % del-shortform
\hline
1 & $1$ & $\aaa \varepsilon^{0} \oplus \aaa \varepsilon^{1} \oplus \chi_{7}^{5} \varepsilon^{2} \oplus \aaa \varepsilon^{8}$ \\
% \hline %% del-shortform
% \hline %% del-shortform
%  & $T_{2}$ & $(1-2x)(1-4x)(1-256x)(1+1294x)$ \\ % del-shortform
%  & $T_{3}$ & $(1-3x)(1-9x)(1+1293x)(1+5476x)$ \\ % del-shortform
%  & $T_{5}$ & $1+5272x + O(x^2)$ \\ % del-shortform
\hline
1 & $1$ & $\chi_{7}^{5} \varepsilon^{0} \oplus \aaa \varepsilon^{1} \oplus \aaa \varepsilon^{2} \oplus \aaa \varepsilon^{8}$ \\
% \hline %% del-shortform
% \hline %% del-shortform
%  & $T_{2}$ & $(1-x)(1-2x)(1-5035x)(1-5123x)$ \\ % del-shortform
%  & $T_{3}$ & $(1-x)(1-3x)(1-4122x)(1-5052x)$ \\ % del-shortform
%  & $T_{5}$ & $1-4947x + O(x^2)$ \\ % del-shortform
\hline
1 & $1$ & $\aaa \varepsilon^{0} \oplus \aaa \varepsilon^{1} \oplus \aaa \varepsilon^{2} \sigma_{7.7.\mathrm{d}.\mathrm{a}}$ \\
\hline
\end{tabular}
\end{center}

\begin{center}
\begin{tabular}{|r|r|c|}
\hline
\multicolumn{3}{|l|}{\textbf{Level $N = 7$.  Coeffs $\Sym^5(V)\otimes\chi_7^3$.  Field $\mathbb{F} = GF(12037^{6})$.}} \\
\multicolumn{3}{|l|}{Computed $T_{2}$, $T_{3}$, $T_{5}$.  Dim~$6$.} \\ % del-longform
% \hline %% del-shortform
%  & $T_{2}$ & $(1-x)(1-2x)(1-4x)(1-256x)$ \\ % del-shortform
%  & $T_{3}$ & $(1-x)(1-3x)(1+9x)(1+5476x)$ \\ % del-shortform
%  & $T_{5}$ & $(1-x)(1-5x)(1+25x)(1-5441x)$ \\ % del-shortform
\hline
1 & $1$ & $\aaa \varepsilon^{0} \oplus \aaa \varepsilon^{1} \oplus \chi_{7}^{3} \varepsilon^{2} \oplus \aaa \varepsilon^{8}$ \\
% \hline %% del-shortform
% \hline %% del-shortform
%  & $T_{2}$ & $(1-x)(1-2x)(1-4x)(1-256x)$ \\ % del-shortform
%  & $T_{3}$ & $(1+x)(1-3x)(1-9x)(1+5476x)$ \\ % del-shortform
%  & $T_{5}$ & $(1+x)(1-5x)(1-25x)(1-5441x)$ \\ % del-shortform
\hline
1 & $1$ & $\chi_{7}^{3} \varepsilon^{0} \oplus \aaa \varepsilon^{1} \oplus \aaa \varepsilon^{2} \oplus \aaa \varepsilon^{8}$ \\
% \hline %% del-shortform
% \hline %% del-shortform
%  & $T_{2}$ & $(1-2x)(1-256x)(1-3557x)(1+3560x)$ \\ % del-shortform
%  & $T_{3}$ & $(1+3x)(1-3x)^2(1+5476x)$ \\ % del-shortform
%  & $T_{5}$ & $(1+5x)(1-5x)^2(1-5441x)$ \\ % del-shortform
\hline
1 & $1$ & $\aaa \varepsilon^{1} \oplus \aaa \varepsilon^{8} \oplus \aaa \varepsilon^{0} \sigma_{7.3.\mathrm{b}.\mathrm{a}}$ \\
\hline
\end{tabular}
\end{center}

\begin{center}
\begin{tabular}{|r|r|c|}
\hline
\multicolumn{3}{|l|}{\textbf{Level $N = 7$.  Coeffs $\Sym^6(V)$.  Field $\mathbb{F} = GF(12037^{6})$.}} \\
\multicolumn{3}{|l|}{Computed $T_{2}$, $T_{3}$, $T_{5}$.  Dim~$3$.} \\ % del-longform
% \hline %% del-shortform
%  & $T_{2}$ & $(1-x)(1-2x)(1-3666x)(1+3690x)$ \\ % del-shortform
%  & $T_{3}$ & $(1-x)(1-3x)(1+378x-3408x^{2})$ \\ % del-shortform
%  & $T_{5}$ & $(1-x)(1-5x)(1-4768x)(1-5169x)$ \\ % del-shortform
\hline
1 & $1$ & $\aaa \varepsilon^{0} \oplus \aaa \varepsilon^{1} \oplus \aaa \varepsilon^{2} \sigma_{7.8.\mathrm{a}.\mathrm{a}}$ \\
%% \hline
%% \multicolumn{3}{|l|}{$\sigma_{7} = \mathrm{Cusp\ form\ }q - 6q^{2} - 42q^{3} - 92q^{4} - 84q^{5} + 252q^{6} + 343q^{7} + O(q^{8})$} \\
%% \multicolumn{3}{|l|}{$\mathrm{for\ }\Gamma_0(7)\mathrm{\ with\ character\ }\aaa\mathrm{\ in\ }\mathbb{F}$} \\
%% \multicolumn{3}{|l|}{$\mathrm{of\ weight\ }8\mathrm{\ with\ coeffs\ in\ }\mathbb{Q}$} \\
% \hline %% del-shortform
% \hline %% del-shortform
%  & $T_{2}$ & $(1 + \left(292 c + 3820\right)x)(1 + \left(8722 c + 11794\right)x)(1 + \left(9680 c + 800\right)x + \left(10536 c + 2559\right)x^{2}) \mathrm{\ with\ } c^2 + 12036*c + 5 = 0$ \\ % del-shortform
%  & $T_{3}$ & $(1 + \left(12036 c + 1961\right)x)(1 + \left(11835 c + 8260\right)x + \left(2755 c + 9313\right)x^{2} + \left(7264 c + 9706\right)x^{3}) \mathrm{\ with\ } c^2 + 12036*c + 5 = 0$ \\ % del-shortform
%  & $T_{5}$ & $(1 + \left(1271 c + 5888\right)x)(1 + \left(1596 c + 6325\right)x)(1 + \left(11559 c + 7393\right)x + \left(7935 c + 3913\right)x^{2}) \mathrm{\ with\ } c^2 + 12036*c + 5 = 0$ \\ % del-shortform
\hline
2 & $1$ & $\aaa \varepsilon^{0} \oplus \aaa \varepsilon^{1} \oplus \aaa \varepsilon^{2} \sigma_{7.8.\mathrm{a}.\mathrm{b}}$ \\
\hline
\end{tabular}
\end{center}

\begin{center}
\begin{tabular}{|r|r|c|}
\hline
\multicolumn{3}{|l|}{\textbf{Level $N = 7$.  Coeffs $\Sym^6(V)\otimes\chi_7^2$.  Field $\mathbb{F} = GF(12037^{6})$.}} \\
\multicolumn{3}{|l|}{Computed $T_{2}$, $T_{3}$, $T_{5}$.  Dim~$4$.} \\ % del-longform
% \hline %% del-shortform
%  & $T_{2}$ & $(1-x)(1-2x)(1+1880x-76x^{2})$ \\ % del-shortform
%  & $T_{3}$ & $(1-x)(1-3x)(1-5345x+4410x^{2})$ \\ % del-shortform
%  & $T_{5}$ & $(1-x)(1-5x)(1-134x+2479x^{2})$ \\ % del-shortform
\hline
4 & $1$ & $\aaa \varepsilon^{0} \oplus \aaa \varepsilon^{1} \oplus \aaa \varepsilon^{2} \sigma_{7.8.\mathrm{c}.\mathrm{a}}$ \\
\hline
\end{tabular}
\end{center}

\begin{center}
\begin{tabular}{|r|r|c|}
\hline
\multicolumn{3}{|l|}{\textbf{Level $N = 7$.  Coeffs $\Sym^7(V)\otimes\chi_7$.  Field $\mathbb{F} = GF(12037^{6})$.}} \\
\multicolumn{3}{|l|}{Computed $T_{2}$, $T_{3}$, $T_{5,1}$.  Dim~$6$.} \\ % del-longform
% \hline %% del-shortform
%  & $T_{2}$ & $(1-x)(1-2x)(1-1024x)(1+5176x)$ \\ % del-shortform
%  & $T_{3}$ & $(1-x)(1-3x)(1-400x)(1+1136x)$ \\ % del-shortform
%  & $T_{5}$ & $1+137x + O(x^2)$ \\ % del-shortform
\hline
1 & $1$ & $\aaa \varepsilon^{0} \oplus \aaa \varepsilon^{1} \oplus \chi_{7} \varepsilon^{2} \oplus \aaa \varepsilon^{10}$ \\
% \hline %% del-shortform
% \hline %% del-shortform
%  & $T_{2}$ & $(1-2x)(1-4x)(1-1024x)(1+1294x)$ \\ % del-shortform
%  & $T_{3}$ & $(1-3x)(1-9x)(1+1136x)(1+1293x)$ \\ % del-shortform
%  & $T_{5}$ & $1-4942x + O(x^2)$ \\ % del-shortform
\hline
1 & $1$ & $\chi_{7} \varepsilon^{0} \oplus \aaa \varepsilon^{1} \oplus \aaa \varepsilon^{2} \oplus \aaa \varepsilon^{10}$ \\
% \hline %% del-shortform
% \hline %% del-shortform
%  & $T_{2}$ & $(1-x)(1-2x)(1-2122x)(1-5716x)$ \\ % del-shortform
%  & $T_{3}$ & $(1-x)(1-3x)(1-4287x+3006x^{2})$ \\ % del-shortform
%  & $T_{5}$ & $1+3795x + O(x^2)$ \\ % del-shortform
\hline
4 & $1$ & $\aaa \varepsilon^{0} \oplus \aaa \varepsilon^{1} \oplus \aaa \varepsilon^{2} \sigma_{7.9.\mathrm{d}.\mathrm{a}}$ \\
\hline
\end{tabular}
\end{center}

\begin{center}
\begin{tabular}{|r|r|c|}
\hline
\multicolumn{3}{|l|}{\textbf{Level $N = 7$.  Coeffs $\Sym^7(V)\otimes\chi_7^3$.  Field $\mathbb{F} = GF(12037^{6})$.}} \\
\multicolumn{3}{|l|}{Computed $T_{2}$, $T_{3}$.  Dim~$8$.} \\ % del-longform
% \hline %% del-shortform
%  & $T_{2}$ & $(1-x)(1-2x)(1-4x)(1-1024x)$ \\ % del-shortform
%  & $T_{3}$ & $(1-x)(1-3x)(1+9x)(1+1136x)$ \\ % del-shortform
\hline
1 & $1$ & $\aaa \varepsilon^{0} \oplus \aaa \varepsilon^{1} \oplus \chi_{7}^{3} \varepsilon^{2} \oplus \aaa \varepsilon^{10}$ \\
% \hline %% del-shortform
% \hline %% del-shortform
%  & $T_{2}$ & $(1-x)(1-2x)(1-4x)(1-1024x)$ \\ % del-shortform
%  & $T_{3}$ & $(1+x)(1-3x)(1-9x)(1+1136x)$ \\ % del-shortform
\hline
1 & $1$ & $\chi_{7}^{3} \varepsilon^{0} \oplus \aaa \varepsilon^{1} \oplus \aaa \varepsilon^{2} \oplus \aaa \varepsilon^{10}$ \\
% \hline %% del-shortform
% \hline %% del-shortform
%  & $T_{2}$ & $(1-2x)(1-1024x)(1-3557x)(1+3560x)$ \\ % del-shortform
%  & $T_{3}$ & $(1+3x)(1-3x)^2(1+1136x)$ \\ % del-shortform
\hline
1 & $1$ & $\aaa \varepsilon^{1} \oplus \aaa \varepsilon^{10} \oplus \aaa \varepsilon^{0} \sigma_{7.3.\mathrm{b}.\mathrm{a}}$ \\
%% \hline
%% \multicolumn{3}{|l|}{$\sigma_{7} = \mathrm{Cusp\ form\ }q - 3q^{2} + 5q^{4} - 7q^{7} + O(q^{8})$} \\
%% \multicolumn{3}{|l|}{$\mathrm{for\ }\Gamma_1(7)\mathrm{\ with\ character\ }\chi_{7}^{3}\mathrm{\ in\ }\mathbb{F}$} \\
%% \multicolumn{3}{|l|}{$\mathrm{of\ weight\ }3\mathrm{\ with\ coeffs\ in\ }\mathbb{Q}$} \\
% \hline %% del-shortform
% \hline %% del-shortform
%  & $T_{2}$ & $(1-x)(1-2x)(1-5384x)(1+5508x)$ \\ % del-shortform
%  & $T_{3}$ & $(1-x)(1-3x)(1+729x)(1-729x)$ \\ % del-shortform
\hline
1 & $1$ & $\aaa \varepsilon^{0} \oplus \aaa \varepsilon^{1} \oplus \aaa \varepsilon^{2} \sigma_{7.9.\mathrm{b}.\mathrm{a}}$ \\
\hline
\end{tabular}
\end{center}

\subsubsection*{Level 8}

\begin{center}
\begin{tabular}{|r|r|c|}
\hline
\multicolumn{3}{|l|}{\textbf{Level $N = 8$.  Coeffs $V\otimes\chi_{8,0}$.  Field $\mathbb{F} = GF(12037^{2})$.}} \\
\multicolumn{3}{|l|}{Computed $T_{3}$, $T_{5}$, $T_{7}$.  Dim~$6$.} \\ % del-longform
% \hline %% del-shortform
%  & $T_{3}$ & $(1-x)(1-3x)(1+9x)(1-81x)$ \\ % del-shortform
%  & $T_{5}$ & $(1-x)(1-5x)(1-25x)(1-625x)$ \\ % del-shortform
%  & $T_{7}$ & $(1-x)(1-7x)(1+49x)(1-2401x)$ \\ % del-shortform
\hline
1 & $3$ & $\aaa \varepsilon^{0} \oplus \aaa \varepsilon^{1} \oplus \chi_{8,0} \varepsilon^{2} \oplus \aaa \varepsilon^{4}$ \\
% \hline %% del-shortform
% \hline %% del-shortform
%  & $T_{3}$ & $(1+x)(1-3x)(1-9x)(1-81x)$ \\ % del-shortform
%  & $T_{5}$ & $(1-x)(1-5x)(1-25x)(1-625x)$ \\ % del-shortform
%  & $T_{7}$ & $(1+x)(1-7x)(1-49x)(1-2401x)$ \\ % del-shortform
\hline
1 & $3$ & $\chi_{8,0} \varepsilon^{0} \oplus \aaa \varepsilon^{1} \oplus \aaa \varepsilon^{2} \oplus \aaa \varepsilon^{4}$ \\
\hline
\end{tabular}
\end{center}

\begin{center}
\begin{tabular}{|r|r|c|}
\hline
\multicolumn{3}{|l|}{\textbf{Level $N = 8$.  Coeffs $V\otimes\chi_{8,0}\chi_{8,1}$.  Field $\mathbb{F} = GF(12037^{2})$.}} \\
\multicolumn{3}{|l|}{Computed $T_{3}$, $T_{5}$, $T_{7}$.  Dim~$4$.} \\ % del-longform
% \hline %% del-shortform
%  & $T_{3}$ & $(1-x)(1-3x)(1-9x)(1-81x)$ \\ % del-shortform
%  & $T_{5}$ & $(1-x)(1-5x)(1+25x)(1-625x)$ \\ % del-shortform
%  & $T_{7}$ & $(1-x)(1-7x)(1+49x)(1-2401x)$ \\ % del-shortform
\hline
1 & $1$ & $\aaa \varepsilon^{0} \oplus \aaa \varepsilon^{1} \oplus \chi_{8,0}\chi_{8,1} \varepsilon^{2} \oplus \aaa \varepsilon^{4}$ \\
% \hline %% del-shortform
% \hline %% del-shortform
%  & $T_{3}$ & $(1-x)(1-3x)(1-9x)(1-81x)$ \\ % del-shortform
%  & $T_{5}$ & $(1+x)(1-5x)(1-25x)(1-625x)$ \\ % del-shortform
%  & $T_{7}$ & $(1+x)(1-7x)(1-49x)(1-2401x)$ \\ % del-shortform
\hline
1 & $1$ & $\chi_{8,0}\chi_{8,1} \varepsilon^{0} \oplus \aaa \varepsilon^{1} \oplus \aaa \varepsilon^{2} \oplus \aaa \varepsilon^{4}$ \\
% \hline %% del-shortform
% \hline %% del-shortform
%  & $T_{3}$ & $(1-x)(1-3x)(1+18x+729x^{2})$ \\ % del-shortform
%  & $T_{5}$ & $(1-x)(1-5x)(1+125x)(1-125x)$ \\ % del-shortform
%  & $T_{7}$ & $(1-x)(1-7x)(1+343x)(1-343x)$ \\ % del-shortform
\hline
1 & $1$ & $\aaa \varepsilon^{0} \oplus \aaa \varepsilon^{1} \oplus \aaa \varepsilon^{2} \sigma_{8.3.\mathrm{d}.\mathrm{a}}$ \\
\hline
\end{tabular}
\end{center}

\begin{center}
\begin{tabular}{|r|r|c|}
\hline
\multicolumn{3}{|l|}{\textbf{Level $N = 8$.  Coeffs $\Sym^2(V)$.  Field $\mathbb{F} = GF(12379^{4})$.}} \\
\multicolumn{3}{|l|}{Computed $T_{3}$, $T_{5}$, $T_{7}$, $T_{11}$.  Dim~$1$.} \\ % del-longform
% \hline %% del-shortform
%  & $T_{3}$ & $(1-x)(1-3x)(1+36x+2187x^{2})$ \\ % del-shortform
%  & $T_{5}$ & $(1-x)(1-5x)(1-2268x)(1+2318x)$ \\ % del-shortform
%  & $T_{7}$ & $(1-x)(1-7x)(1-1176x-5850x^{2})$ \\ % del-shortform
%  & $T_{11}$ & $(1-x)(1-11x)(1+5324x+2625x^{2})$ \\ % del-shortform
\hline
1 & $1$ & $\aaa \varepsilon^{0} \oplus \aaa \varepsilon^{1} \oplus \aaa \varepsilon^{2} \sigma_{8.4.\mathrm{a}.\mathrm{a}}$ \\
%% \hline
%% \multicolumn{3}{|l|}{$\sigma_{8} = \mathrm{Cusp\ form\ }q - 4q^{3} - 2q^{5} + 24q^{7} - 11q^{9} - 44q^{11} + O(q^{12})$} \\
%% \multicolumn{3}{|l|}{$\mathrm{for\ }\Gamma_0(8)\mathrm{\ with\ character\ }\aaa\mathrm{\ in\ }\mathbb{F}$} \\
%% \multicolumn{3}{|l|}{$\mathrm{of\ weight\ }4\mathrm{\ with\ coeffs\ in\ }\mathbb{Q}$} \\
\hline
\end{tabular}
\end{center}

\begin{center}
\begin{tabular}{|r|r|c|}
\hline
\multicolumn{3}{|l|}{\textbf{Level $N = 8$.  Coeffs $\Sym^2(V)\otimes\chi_{8,1}$.  Field $\mathbb{F} = GF(12037^{2})$.}} \\
\multicolumn{3}{|l|}{Computed $T_{3}$, $T_{5}$, $T_{7}$.  Dim~$2$.} \\ % del-longform
% \hline %% del-shortform
%  & $T_{3}$ & $(1-x)(1-3x)(1-4301x-2187x^{2})$ \\ % del-shortform
%  & $T_{5}$ & $(1-x)(1-5x)(1+4461x)(1-5978x)$ \\ % del-shortform
%  & $T_{7}$ & $(1-x)(1-7x)(1+392x+5027x^{2})$ \\ % del-shortform
\hline
2 & $1$ & $\aaa \varepsilon^{0} \oplus \aaa \varepsilon^{1} \oplus \aaa \varepsilon^{2} \sigma_{8.4.\mathrm{b}.\mathrm{a}}$ \\
\hline
\end{tabular}
\end{center}

\subsubsection*{Level 9}

\begin{center}
\begin{tabular}{|r|r|c|}
\hline
\multicolumn{3}{|l|}{\textbf{Level $N = 9$.  Coeffs $V\otimes\chi_9$.  Field $\mathbb{F} = GF(12037^{2})$.}} \\
\multicolumn{3}{|l|}{Computed $T_{2}$, $T_{5}$, $T_{7}$.  Dim~$4$.} \\ % del-longform
% \hline %% del-shortform
%  & $T_{2}$ & $(1-x)(1-2x)(1-16x)(1+5172x)$ \\ % del-shortform
%  & $T_{5}$ & $(1-x)(1-5x)(1-625x)(1+3761x)$ \\ % del-shortform
%  & $T_{7}$ & $(1-x)(1-7x)(1-2401x)(1-3172x)$ \\ % del-shortform
\hline
1 & $1$ & $\aaa \varepsilon^{0} \oplus \aaa \varepsilon^{1} \oplus \chi_{9} \varepsilon^{2} \oplus \aaa \varepsilon^{4}$ \\
% \hline %% del-shortform
% \hline %% del-shortform
%  & $T_{2}$ & $(1-2x)(1-4x)(1-16x)(1+1293x)$ \\ % del-shortform
%  & $T_{5}$ & $(1-5x)(1-25x)(1-625x)(1-1294x)$ \\ % del-shortform
%  & $T_{7}$ & $(1-7x)(1-49x)(1-1293x)(1-2401x)$ \\ % del-shortform
\hline
1 & $1$ & $\chi_{9} \varepsilon^{0} \oplus \aaa \varepsilon^{1} \oplus \aaa \varepsilon^{2} \oplus \aaa \varepsilon^{4}$ \\
% \hline %% del-shortform
% \hline %% del-shortform
%  & $T_{2}$ & $(1-x)(1-2x)(1+223x)(1-5391x)$ \\ % del-shortform
%  & $T_{5}$ & $(1-x)(1-5x)(1-4565x-3410x^{2})$ \\ % del-shortform
%  & $T_{7}$ & $(1-x)(1-7x)(1+2766x)(1+3676x)$ \\ % del-shortform
\hline
1 & $1$ & $\aaa \varepsilon^{0} \oplus \aaa \varepsilon^{1} \oplus \aaa \varepsilon^{2} \sigma_{9.3.\mathrm{d}.\mathrm{a}}$ \\
\hline
\end{tabular}%
\end{center}

\subsubsection*{Level 10}

\begin{center}
%
\end{center}

\subsubsection*{Level 11}

\begin{center}
%
\end{center}

\subsubsection*{Level 12}

\begin{center}
%
\end{center}

\subsubsection*{Level 13}

\begin{center}
%
\end{center}

\subsubsection*{Level 14}

\begin{center}
%
\end{center}

\subsubsection*{Level 15}

\begin{center}
%
\end{center}

\subsubsection*{Level 16}

\begin{center}
%
\end{center}

\subsubsection*{Level 17}

\begin{center}
%
\end{center}

\subsubsection*{Level 18}

\begin{center}
%
\end{center}

\subsection{} \label{resChi}
For each~$N$, the next table specifies the basis that Sage
chooses for the group of characters $(\Z/N\Z)^\times \to \F_p$.  If
there is one basis element, it is denoted~$\chi_N$.  If there is more
than one, they are denoted $\chi_{N,0}$, $\chi_{N,1}$, etc.  The
\emph{order} of~$\chi$ is the smallest positive~$n$ so that $\chi^n$
is trivial on $(\Z/N\Z)^\times$.  The \emph{parity} is even if
$\chi(-1) = +1$ and odd if $\chi(-1) = -1$.

\begin{center}
%
\end{center}

\bibliographystyle{amsalpha_no_mr}
\bibliography{AGM-VIII}

\def\Dbar{\leavevmode\lower.6ex\hbox to 0pt{\hskip-.23ex \accent"16\hss}D}
  \def\cftil#1{\ifmmode\setbox7\hbox{$\accent"5E#1$}\else
  \setbox7\hbox{\accent"5E#1}\penalty 10000\relax\fi\raise 1\ht7
  \hbox{\lower1.15ex\hbox to 1\wd7{\hss\accent"7E\hss}}\penalty 10000
  \hskip-1\wd7\penalty 10000\box7}
  \def\cfudot#1{\ifmmode\setbox7\hbox{$\accent"5E#1$}\else
  \setbox7\hbox{\accent"5E#1}\penalty 10000\relax\fi\raise 1\ht7
  \hbox{\raise.1ex\hbox to 1\wd7{\hss.\hss}}\penalty 10000 \hskip-1\wd7\penalty
  10000\box7}
\providecommand{\bysame}{\leavevmode\hbox to3em{\hrulefill}\thinspace}
\providecommand{\MR}{\relax\ifhmode\unskip\space\fi MR }
\providecommand{\MRhref}[2]{%
  \href{http://www.ams.org/mathscinet-getitem?mr=#1}{#2}
}
\providecommand{\href}[2]{#2}
\begin{thebibliography}{{LMF}23}

\bibitem[ADP02]{ADP}
Avner Ash, Darrin Doud, and David Pollack, \emph{Galois representations with
  conjectural connections to arithmetic cohomology}, Duke Math. J. \textbf{112}
  (2002), no.~3, 521--579.

\bibitem[AGM02]{AGM1}
Avner Ash, Paul~E. Gunnells, and Mark McConnell, \emph{Cohomology of congruence
  subgroups of {${\rm SL}_4(\mathbb{Z})$}}, J. Number Theory \textbf{94}
  (2002), no.~1, 181--212.

\bibitem[AGM08]{AGM2}
\bysame, \emph{Cohomology of congruence subgroups of {${\rm
  SL}_4(\mathbb{Z})$}. {II}}, J. Number Theory \textbf{128} (2008), no.~8,
  2263--2274.

\bibitem[AGM10]{AGM3}
\bysame, \emph{Cohomology of congruence subgroups of {${\rm
  SL}_4(\mathbb{Z})$}. {III}}, Math. Comp. \textbf{79} (2010), no.~271,
  1811--1831.

\bibitem[AGM11]{AGM4}
\bysame, \emph{Torsion in the cohomology of congruence subgroups of
  $\mathrm{SL}(4,\mathbb{Z})$ and {G}alois representations}, J. Algebra
  \textbf{325} (2011), 404--415.

\bibitem[AGM20]{AGM7}
\bysame, \emph{Cohomology with twisted one-dimensional coefficients for
  congruence subgroups of {${\rm SL}_4(\Bbb Z)$} and {G}alois representations},
  J. Algebra \textbf{553} (2020), 211--247.

\bibitem[CLH16]{CLH}
Ana Caraiani and Bao~V. Le~Hung, \emph{On the image of complex conjugation in
  certain galois representations}, Compositio Math. \textbf{152} (2016), no.~7,
  1476–1488.

\bibitem[{LMF}23]{lmfdb}
The {LMFDB Collaboration}, \emph{The {L}-functions and modular forms database},
  \texttt{https://www.lmfdb.org}, 2023, [Online; accessed 30 July 2023].

\bibitem[Sch09]{AS}
Achill Sch\"{u}rmann, \emph{Enumerating perfect forms}, Quadratic
  forms---algebra, arithmetic, and geometry, Contemp. Math., vol. 493, Amer.
  Math. Soc., Providence, RI, 2009, pp.~359--377.

\end{thebibliography}

\end{document}